\documentclass[journal]{IEEEtran}
\pdfoutput=1

\ifCLASSINFOpdf
\else
\fi
\usepackage[cmex10]{amsmath}

\usepackage{mathtools} 
\usepackage{amssymb}
\usepackage{multirow}
\usepackage{siunitx}
\usepackage{balance} 

\usepackage[amsmath,amsthm,thmmarks]{ntheorem}
\newtheorem{thmresult}{Result}

\newtheorem{corollary}{Corollary}[thmresult]

\usepackage{float}
\floatstyle{ruled}
\newfloat{algorithm}{tbp}{loa}
\providecommand{\algorithmname}{Algorithm}
\floatname{algorithm}{\protect\algorithmname}

\usepackage[breaklinks]{hyperref}
\hypersetup{
	bookmarks=true,         
	colorlinks=true,       
    linkcolor=blue,          
    citecolor=green,        
    filecolor=magenta,      
    urlcolor=magenta,           
  pdfauthor={Daniel Berjon, Guillermo Gallego, Carlos Cuevas, Francisco Moran, Narciso Garcia},%
  pdftitle={Optimal Piecewise Linear Function Approximation for GPU-based Applications},%
  pdfsubject={Cybernetics, Image processing, Computer Vision},%
  pdfkeywords={Numerical approximation and analysis, error equalization, orthogonal projection, parallel processing, piecewise linearization, Gaussian function, Lorentzian function, Bessel function.},%
  pdfcreationdate={D:20150923120000}
}

\begin{document}

\global\long\def\Ltwo{L_{2}}
\global\long\def\v{v}
\global\long\def\w{w}
\global\long\def\f{f}
\global\long\def\part{T}
\global\long\def\x{\mathbf{x}}
\global\long\def\R{\mathbb{R}}
\global\long\def\vecspace{\mathrm{V}_{\part}}
\global\long\def\pip{\pi_{\part}}
\global\long\def\ortproj{P_{\part}}
\global\long\def\inner#1#2{\left\langle #1,#2\right\rangle }
\global\long\def\dif{\mathrm{d}}
\global\long\def\partOpt{\part^{\ast}}
\global\long\def\partQuasiOpt{\part^{(\ast)}}
\global\long\def\partUnif{\part_U}
\global\long\def\vecspaceOpt{\mathrm{V}_{\partOpt}}
\global\long\def\pipOpt{\pi_{\partOpt}}
\global\long\def\ortprojOpt{P_{\partOpt}}
\global\long\def\pipQuasiOpt{\pi_{\partQuasiOpt}}
\global\long\def\ortprojQuasiOpt{P_{\partQuasiOpt}}
\global\long\def\pipUnif{\pi_{\partUnif}}
\global\long\def\ortprojUnif{P_{\partUnif}}
\global\long\def\xsamp{\mathtt{x\_samples}}
\global\long\def\vsamp{\mathtt{v\_samples}}
\global\long\def\fsec{\f^{\prime\prime}}
\global\long\def\fimax{|\fsec_{i}|_{\max}}
\global\long\def\fZero{\fsec_{0x}}
\global\long\def\fOne{\fsec_{1x}}
\global\long\def\bb{\mathbf{b}}
\global\long\def\by{\mathbf{y}}
\global\long\def\mA{\mathtt{A}}
\global\long\def\leftDeltaY#1{\Delta y_{#1}}
\global\long\def\rightDeltaY#1{\Delta y_{#1}}

\title{Optimal Piecewise Linear Function Approximation for GPU-based Applications}

\author{Daniel~Berj\'{o}n, Guillermo~Gallego, Carlos~Cuevas, Francisco~Mor\'{a}n and~Narciso~Garc\'{i}a
\thanks{Manuscript received March 7, 2014; revised June 26, 2015; accepted
September 22, 2015.
This work has been supported in part by the Ministerio de Econom\'{i}a y Competitividad of the Spanish Government under grant TEC2013-48453 (MR-UHDTV) and by the European Commission under grant 610691 (BRIDGET).}
\thanks{D. Berj\'{o}n, C. Cuevas, F. Mor\'{a}n, and N. Garc\'{i}a are with the Grupo de Tratamiento de Im\'{a}genes, ETSI Telecomunicaci\'{o}n, Universidad Polit\'{e}cnica de Madrid, Madrid 28040, Spain (e-mail: dbd@gti.ssr.upm.es, ccr@gti.ssr.upm.es, fmb@gti.ssr.upm.es, narciso@gti.ssr.upm.es).}
\thanks{G. Gallego was also with the Grupo de Tratamiento de Im\'{a}genes, Universidad Polit\'{e}cnica de Madrid, Madrid 28040, Spain. 
He is now with the Robotics and Perception Group, University of Zurich, Zurich 8001, Switzerland (e-mail: guillermo.gallego@ifi.uzh.ch).

Copyright (c) 2015 IEEE. Personal use of this material is permitted.
However, permission to use this material for any other purposes must
be obtained from the IEEE. 
\protect\href{http://dx.doi.org/10.1109/TCYB.2015.2482365}{DOI: 10.1109/TCYB.2015.2482365}}
}

\maketitle

\begin{abstract}
Many computer vision and human-computer interaction applications developed
in recent years need evaluating complex and continuous mathematical
functions as an essential step toward proper operation. 
However, rigorous evaluation of this kind of functions often implies a very high
computational cost, 
unacceptable in real-time applications. 
To alleviate this problem, functions are commonly approximated by 
simpler piecewise-polynomial representations.
Following this idea, we propose a novel, efficient, and practical 
technique to evaluate complex and continuous functions using
a nearly optimal design of two types of piecewise linear approximations 
in the case of a large budget of evaluation subintervals.
To this end, we develop a thorough 
error analysis that yields asymptotically tight bounds to accurately quantify the 
approximation performance of both representations.
It provides an improvement upon previous error estimates and 
allows the user to control the trade-off between 
the approximation error and the number of evaluation subintervals.
To guarantee real-time operation, the method is suitable for, but
not limited to, an efficient implementation in modern Graphics
Processing Units (GPUs), where it outperforms previous alternative
approaches by exploiting the fixed-function interpolation routines
present in their texture units.
The proposed technique is a perfect match for
any application requiring the evaluation of continuous functions; 
we have measured in detail its quality and efficiency on several functions, and, in particular, the Gaussian 
function because it is extensively used in many areas of computer vision and cybernetics, 
and it is expensive to evaluate.
\end{abstract}


\IEEEpeerreviewmaketitle

\section{Introduction}
\IEEEPARstart{T}{he design} of high-quality and real-time image processing algorithms 
is a key topic in today's context, where humans
demand more capabilities and control over their rapidly increasing 
number of advanced imaging devices such as cameras, smart phones,
tablets, etc.~\cite{Shapiro2013}.  
Many efforts are being made toward satisfying those needs 
by exploiting dedicated hardware, such as Graphics Processing Units
(GPUs)~\cite{SinghalICIP2010,Borges2011,Kapinchev2013}. 
%
Some areas where GPUs are emerging to improve Human Computer Interaction (HCI) include 
cybernetics (for example, in facial~\cite{Song2010} or object~\cite{Kim2009} recognition, 
classification using Support Vector Machines~\cite{HerreroLopez2011} 
and genetic algorithms for clustering~\cite{Kromer2012}), 
and image processing (e.g., unsupervised image segmentation~\cite{Borges2011,BallaArabe2013},
optical coherence tomography systems~\cite{Kapinchev2013},
efficient surface reconstruction from noisy data~\cite{Jalba2009},
remote sensing~\cite{BallaArabe2014},
real-time background subtraction~\cite{Cheng2011}, etc.). 
In such hardware-oriented designed algorithms, 
the computational efficiency of processing tasks is
significantly improved by parallelizing the operations.


\begin{figure}
\begin{minipage}[t]{0.5\columnwidth}%
\begin{center}
\includegraphics[width=0.996\columnwidth,height=0.747\columnwidth]{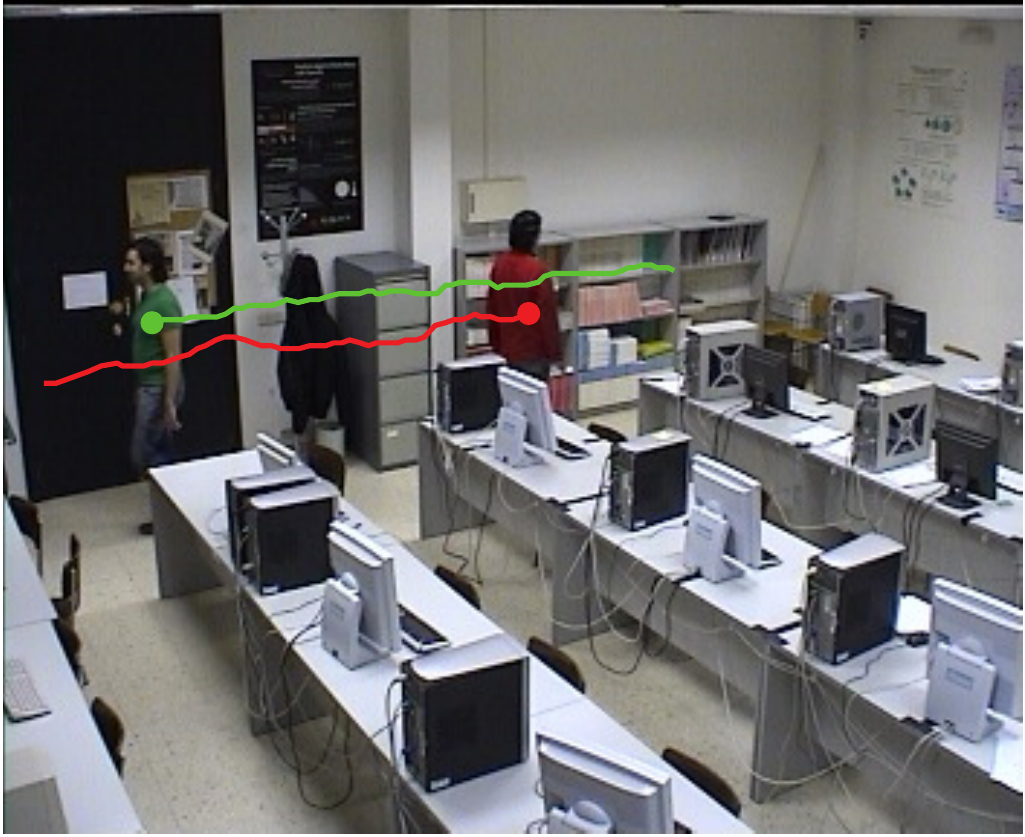}\\
{\footnotesize{} (a)}
\par\end{center}%
\end{minipage}%
\begin{minipage}[t]{0.5\columnwidth}%
\begin{center}
\includegraphics[width=0.996\columnwidth]{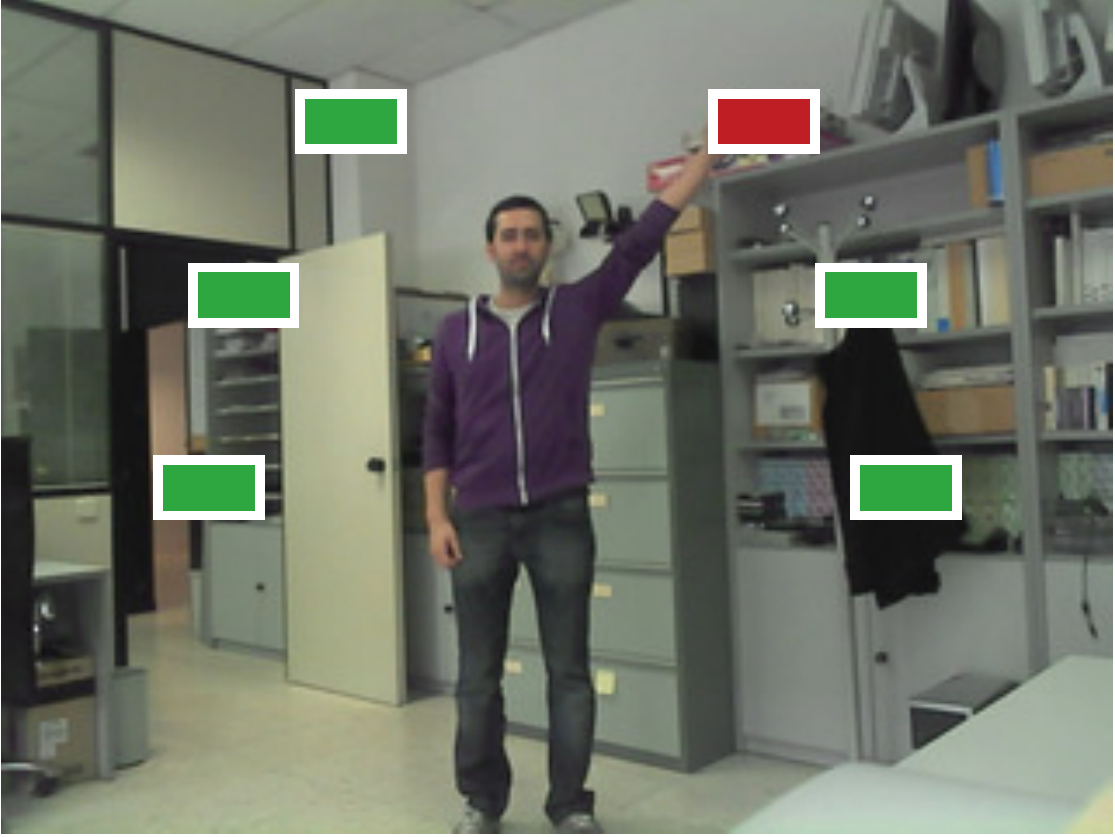}\\
{\footnotesize{} (b)}
\par\end{center}%
\end{minipage}

\caption{\label{fig:applications}High level tasks such as (a) object/person-tracking 
or (b) video-based gestural interfaces often require a foreground segmentation method
as a base building block.}
\end{figure}

Automated video analysis applications such as person tracking or video-based gestural interfaces (see \figurename~\ref{fig:applications}) rely on lower-level building blocks like foreground segmentation~\cite{Chiranjeevi2012robust,Liu2015},
where the design of efficient computational methods for
dedicated hardware has a significant impact.
Multimodal nonparametric segmentation strategies have drawn a lot of
attention~\cite{Sheikh2009background,sigal2012loose} since they are
able to provide high-quality results even in complex scenarios
(dynamic background, illumination changes, 
etc.)~\cite{Cuevas2012}. However, their main
drawback is their extremely high computational cost (requiring the
evaluation of billions of multidimensional Gaussian kernels per
second), which makes them difficult to integrate in the latest
generation of image processing applications~\cite{moshe2012foreground,Cuevas2013}.  
To overcome this important drawback and achieve real-time
performance~\cite{Berjon2013}, the use of parallel hardware such as GPUs
helps but may not be enough by itself, depending on the required resolution,
hence the need for algorithms capable of evaluating non-linear (e.g., Gaussian) functions at high speed
within a required error tolerance.
%
GPU vendors are aware of this recurrent computing problem and
provide hardware implementations of common transcendental functions in
Special Function Units (SFUs)~\cite{Lindholm2008}; indeed, we can find
in the literature examples of successful use of these hardware
facilities~\cite{Sylwestrzak2013}.  However, their ease of use comes
at the price of non-customizable reduced numerical
precision~\cite{NVIDIACorporation2012}.

\subsection{Contribution}
We propose a novel, fast, and practical method to evaluate any 
continuous mathematical function within a known interval. 
Our contributions include, based on error equalization, a nearly optimal design of two 
types of piecewise linear approximations (linear interpolant and orthogonal projection) 
in the $\Ltwo$ norm under the constraint of a large budget of evaluation subintervals~$N$.
Moreover, we provide asymptotically tight bounds for the approximation errors 
of both piecewise linear representations, improving upon existing ones. 
Specifically, in addition to the $O(N^{-2})$ convergence rate 
typical of these approximations, we quantify their interpolation constants:
$1/\sqrt{120}$ for the linear interpolant and 
a further $\sqrt{6}$ improvement factor in case of the orthogonal projection.
The obtained error relations are parameterized by 
the number of segments used to represent the complex (nonlinear) function, 
hence our approach allows the user to estimate the errors given $N$ or, conversely, 
estimate the required $N$ to achieve a target approximation error.

We also propose an efficient implementation of the technique
in a modern GPU by exploiting 
the fixed-function interpolation routines present in its 
texture units to accelerate the computation while 
leaving the rest of the GPU (and, of course, the CPU) free to perform other tasks; 
this technique is even faster than using
SFUs~\cite{Lindholm2008}.

Although the initial motivation of the developed method was to improve
the efficiency of nonparametric foreground segmentation strategies, it
must be noted that it can be also used in many other scientific fields
such as computer
vision~\cite{gallego2008segmentation,guillaumin2012face} or audio
signal processing~\cite{Xie2012,Sehili2012}, where the evaluation of
continuous mathematical functions constitutes a significant computational burden.

\subsection{Organization}
Section~\ref{sec:Piecewise-linear-approximations} summarizes the basic
facts about piecewise linear approximation of real-valued
univariate functions and reviews related work in this topic; 
Section~\ref{sec:Finding-a-non-uniform} derives
a suboptimal partition of the domain of the approximating
function in order to minimize the distance to the original function (proofs of results are given in the Appendixes);
Section~\ref{sec:Computational-analysis-and} analyzes the algorithmic
complexity of the proposed approximation strategies and
Section~\ref{sec:Implementation-in-a-GPU} gives details about their
implementation on modern GPUs. Section~\ref{sec:Results} presents the
experimental results of the proposed approximation on several functions (Gaussian, Lorentzian and Bessel's), both
in terms of quality and computational times, 
and its use is demonstrated in an image processing application. 
Finally, Section~\ref{sec:Conclusions} concludes the paper.

\section{\label{sec:Piecewise-linear-approximations}Piecewise Linear Modelization}

\subsection{Related Work}

In many applications, rigorous evaluation of complex mathematical
functions is not practical because it takes too much computational
power. Consequently, this evaluation is carried out approximating them
by simpler functions such as (piecewise-)polynomial ones.  
Piecewise linearization has been used as an attractive simplified representation of various complex nonlinear systems~\cite{Storace2004}. 
The resulting models fit into well established tools for
linear systems and reduce the complexity of finding the inverse of
nonlinear functions~\cite{Hatanaka2002,Tanjad2011}.  They can also be
used to obtain approximate solutions in complex nonlinear systems, for
example, in Mixed-Integer Linear Programming (MILP)
models~\cite{Rossi2013}.  Some efforts have been devoted as well to
the search for canonical representations in one and multiple
dimensions~\cite{Julian1999,Julian2000} with different goals such as
black box system identification, approximation or model reduction.

Previous attempts to address the problem considered here (optimal function piecewise
linearization) include~\cite{Bellman1969,Frenzen2010437,Hatanaka2002,Ghosh2011}, 
the latter two in the context of nonlinear dynamical systems. 
In~\cite{Bellman1969} an iterative multi-stage procedure
based on dynamical programming is given to provide a solution to the
problem on sequences of progressively finer 2-D
grids. In~\cite{Hatanaka2002} the piecewise linear approximation is
obtained by using the evolutionary computation approach such as
genetic algorithm and evolution strategies; the resulting model is
obtained by minimization of a sampled version of the mean squared
error and it may not be continuous. In~\cite{Ghosh2011} the problem
is addressed using a hybrid approach based on curve fitting,
clustering and genetic algorithms. Here we address the problem from a
less heuristic point of view, using differential calculus to derive a more principled approach~\cite{Gallego2014}.
Our method solely relies on standard numerical integration techniques, 
which takes few seconds to compute, as opposed to recursive partitioning techniques such as~\cite{Frenzen2010437}, which take significantly longer.

\subsection{Basic Results in Piecewise-Linear Interpolation and Least-Squares Approximation}

\begin{figure}
\centering\includegraphics[width=84mm]{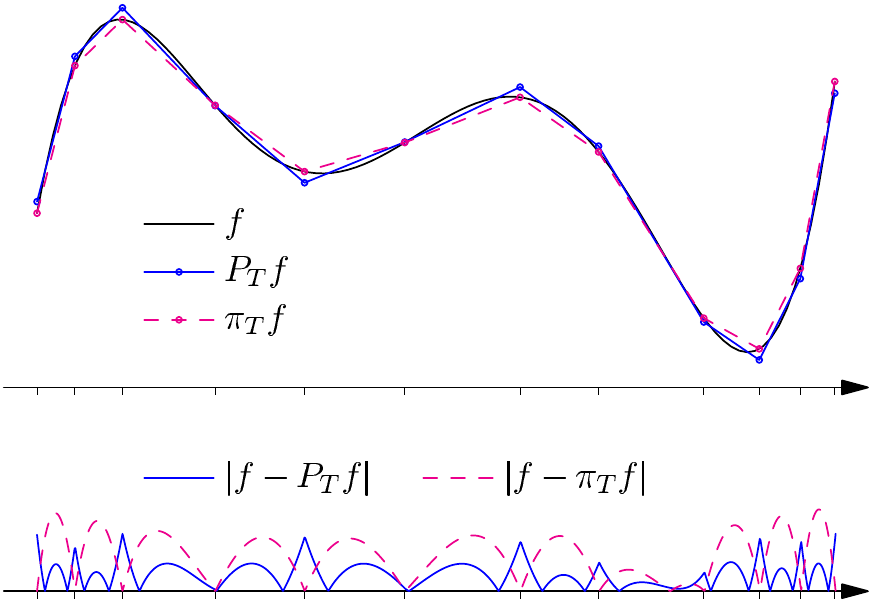}
\caption{\label{fig:pipvsproj}Top: fifth-degree polynomial  
$\f(x)=(x+4)(x+2)(x+1)(x-1)(x-3)$ and two continuous piecewise linear~(CPWL) approximations, the
orthogonal projection $\ortproj\f$ and the linear interpolant
$\pip\f$; bottom: corresponding absolute approximation errors
(magnified by a $5\times$ factor).}
\end{figure}

In this section we summarize the basic theory behind piecewise
linear functions~\cite{Eriksson2004} and two such approximations to
real-valued functions: interpolation and projection
(Sections~\ref{sub:Linear-interpolation}
and~\ref{sub:The-L2-projection}, respectively), 
pictured in \figurename~\ref{fig:pipvsproj} along with their absolute 
approximation error with respect to $\f$.

\begin{figure}
\centering\includegraphics[width=84mm]{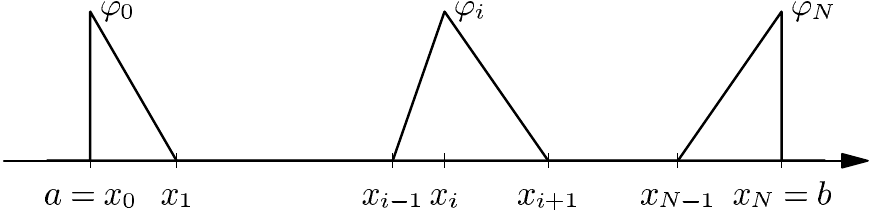}
\caption{\label{fig:Hat-functions-constitute}Hat functions constitute a basis
of $\vecspace$. Since all the basis functions are zero outside
$I=[a,b]$, the basis functions associated with boundary nodes are
only \emph{half hats}.}
\end{figure}

In essence, a piecewise function over an interval $I=\left[a,b\right]$
is a partition of $I$ into a set of $N$ subintervals
$\part=\{I_{i}\}_{i=1}^{N}$, where
$I_{i}=\left(x_{i-1},x_{i}\right)\mid a=x_{0}<x_{1}<\ldots<x_{N}=b$,
and a set of $N$ functions $\f_{i}(x)$, one for each subinterval~$I_{i}$. In particular we are interested in continuous piecewise
linear~(CPWL) functions, which means that all the $\f_{i}(x)$ are
linear and $\f_{i}(x_{i})=\f_{i+1}(x_{i})\;\forall i=\left\{
1,\ldots,N-1\right\} $.
CPWL functions of a given partition $\part$
are elements of a vector space $\vecspace$: the addition of such
functions and/or multiplication by a scalar yields another CPWL
function defined over the same subintervals. A useful basis for the
vector space $\vecspace$ is formed by the set of \emph{hat functions}
or \emph{nodal basis functions} $\left\{ \varphi_{i}\right\}
_{i=0}^{N}$, pictured in \figurename~\ref{fig:Hat-functions-constitute} and
defined in general by the formula
\begin{equation}
\varphi_{i}(x)=\begin{cases}
(x-x_{i-1})/(x_{i}-x_{i-1}), & x\in\left[x_{i-1},x_{i}\right]\\
(x-x_{i+1})/(x_{i}-x_{i+1}), & x\in\left[x_{i},x_{i+1}\right]\\
0, & x\notin\left[x_{i-1},x_{i+1}\right].
\end{cases}
\end{equation}
The basis functions $\varphi_{0}$ and $\varphi_{N}$ associated to the
boundary nodes $x_{0}$ and $x_{N}$ are only \emph{half hats}.

These basis functions are convenient since they can represent any
function $\v$ in $\vecspace$ by just requiring the values of $\v$ at
its nodal points, $\v(x_{i})$, in the form
\begin{equation}
\v(x)=\sum_{i=0}^{N}\v(x_{i})\varphi_{i}(x).
\end{equation}

\subsubsection{\label{sub:Linear-interpolation}Linear Interpolation}

The piecewise linear interpolant $\pip\f\in\vecspace$ of a continuous
function $\f$ over the interval $I$ can be defined in terms of the
basis just introduced:
\begin{equation}
\pip\f(x)=\sum_{i=0}^{N}\f(x_{i})\varphi_{i}(x).
\end{equation}
While this CPWL approximation is trivial to construct, and
may be suitable for some uses, it is by no means the best
possible one. Crucially, $\pip\f(x)\geq\f(x)\;\forall x\in I$ for any
function $\f$ that is convex in $I$. Depending on the application in
which this approximation is used, this property could skew the results. 
However, as we will see in
Section~\ref{sec:Finding-a-non-uniform}, the linear interpolant is
useful to analyze other possible approximations.  It is also at the
heart of the trapezoidal rule for numerical integration.

\subsubsection{\label{sub:The-L2-projection}Orthogonal Projection onto Vector Space $\vecspace$}

Let the usual inner product between two square-integrable ($\Ltwo$)
functions in the interval $I$ be given by
\begin{equation}
\inner uv=\int_{I}u(x)\, v(x)\,\dif x.\label{eq:innerprodDef}
\end{equation}
 Then, the vector space $\vecspace$ can be endowed with the above
inner product, yielding an inner product space. As usual, let
$\|u\|=\sqrt{\inner uu}$ be the norm induced by the inner
product, and let $d(u,v)=\|u-v\|$ be the distance between two
functions $u,v$. The orthogonal projection of the function $\f$ onto
$\vecspace$ is the function $\ortproj\f\in\vecspace$ such that
\begin{equation}
\inner{\f-\ortproj\f}{\v}=0\quad\forall\v\in\vecspace.\label{eq:phf}
\end{equation}
Since $\ortproj\f\in\vecspace$, it can be expressed in the nodal
basis $\left\{ \varphi_{i}\right\} _{i=0}^{N}$ by
\begin{equation}
\ortproj\f(x)=\sum_{i=0}^{N}c_{i}\varphi_{i}(x),\label{eq:OrtProjRepresentation}
\end{equation}
where the coefficients $c_{i}$ solve the linear system of equations
$Mc=b$. The $\left(N+1\right)\times(N+1)$ matrix $M=(m_{ij})$ has
entries $m_{ij}=\inner{\varphi_{i}}{\varphi_{j}}$, and vector $b$ has
entries given by $b_{i}=\inner{\f}{\varphi_{i}}$. The Gramian matrix
$M$ is tridiagonal and strictly diagonally dominant. Therefore, the
system has exactly one solution $c$, which can be obtained efficiently
using the Thomas algorithm (a particular case of Gaussian
elimination)~\cite{deBoor2001}.

$\ortproj\f$ is the element in $\vecspace$ that is closest to $\f$ in
the sense given by the aforementioned $\Ltwo$ distance $d$, as we recall next. 
For any $\w\in\vecspace$,
\[
\begin{aligned}
\|\f-\ortproj\f\|^{2} 
&=\inner{\f-\ortproj\f}{\f-w+w-\ortproj\f}\\
&=\inner{\f-\ortproj\f}{\f-\w}+\inner{\f-\ortproj\f}{\w-\ortproj\f}\\
&\stackrel{\eqref{eq:phf}}{=}\inner{\f-\ortproj\f}{\f-\w},
\end{aligned}
\]
where we used property~\eqref{eq:phf} that the vector
$(\w-\ortproj\f)\in\vecspace$.  Next, applying the Cauchy-Schwarz
inequality,
\[
\|\f-\ortproj\f\|^{2}=\inner{\f-\ortproj\f}{\f-\w}\leq\|\f-\ortproj\f\|\,\|\f-\w\|,
\]
and so $\|\f-\ortproj\f\|\leq\|\f-\w\|$, i.e.,
\begin{equation}
d(\f,\ortproj\f)\leq
d(\f,\w)\quad\forall\w\in\vecspace,\label{eq:phf-is-the-best}
\end{equation}
with equality if $\w=\ortproj\f$. This makes $\ortproj\f$ most
suitable as an approximation of $\f$ under the $\Ltwo$ norm.

\section{\label{sec:Finding-a-non-uniform}Finding the Optimal Partition}

As just established, for a given vector space of CPWL functions
$\vecspace$, $\ortproj\f$ is the function in $\vecspace$ whose $\Ltwo$
distance to $\f$ is minimal. However, the approximation error
$\|\f-\ortproj\f\|$ does not take the same value for every possible
$\vecspace$; therefore, we would like to find the optimal partition
$\partOpt$ (or a reasonable approximation thereof) corresponding to
the space $\vecspaceOpt$ in which the approximation error is minimum,
\[
\|\f-\ortprojOpt\f\|\leq\|\f-\ortproj\f\|.
\]

This is a difficult optimization problem: in order to properly specify
$\ortproj\f$ in \eqref{eq:OrtProjRepresentation} and measure its
distance to $\f$, we need to solve the linear system of equations
$Mc=b$, whose coefficients depend on the partition itself, which makes
the problem symbolically intractable. Numerical algorithms could be
used but it is still a challenging problem.

Let us examine the analogous problem for the interpolant function
$\pip\f(x)$ defined in Section~\ref{sub:Linear-interpolation}.  Again,
we would like to find the optimal partition $\partOpt$ corresponding
to the space $\vecspaceOpt$ in which the approximation error is
minimum,
\[
\|\f-\pipOpt\f\|\leq\|\f-\pip\f\|.
\]
Albeit not as difficult as the previous problem, because $\pip\f$ is
more straightforward to define, this is still a challenging non-linear
optimization problem. Fortunately, as it is shown next, a good
approximation can be easily found under some asymptotic analysis.

In the rest of this section we investigate in detail the error incurred when
approximating a function $\f$ by the interpolant $\pip\f$ and the
orthogonal projection $\ortproj\f$ defined in
Section~\ref{sec:Piecewise-linear-approximations} (see \figurename~\ref{fig:pipvsproj}).  
Then we derive an approximation to the optimal partition (Section~\ref{sub:Approximation-to-the-optimal-partition}) that serves equally well for
$\pip\f$ and $\ortproj\f$ because, as we will show (Section~\ref{sub:Error-in-a-Single-Interval}), their
approximation errors are roughly proportional under the assumption of
a sufficiently large number of intervals.

\subsection{Error in a Single Interval}
\label{sub:Error-in-a-Single-Interval}

\subsubsection{Linear Interpolant}

First, let us consider the error incurred when approximating a
function $\f$, twice continuously differentiable, by its linear
interpolant in an interval (\figurename~\ref{fig:errorzone}).
\begin{thmresult}\label{thm:LinInterpolantErrBound} 
The $\Ltwo$ error between a given function $f$ and its linear interpolant 
\begin{equation}
\pip\f(x)=\f(x_{i-1})\bigl(1-\delta_{i}(x)\bigr)+\f(x_{i})\delta_{i}(x),\label{eq:LinearInterpolantInInterval}
\end{equation}
with $\delta_{i}(x)=(x-x_{i-1})/h_{i}$, 
in the interval $I_{i}=[x_{i-1},x_{i}]$, of length $h_{i}=x_{i}-x_{i-1}$, 
is bounded:
\begin{equation}
\|\f-\pip\f\|_{\Ltwo(I_{i})}\leq\frac{1}{\sqrt{120}}\fimax h_{i}^{5/2},\label{eq:interpolation-bound-for-error-in-an-interval}
\end{equation}
where $\fimax=\max_{\eta\in I_{i}}|\fsec(\eta)|$.
\end{thmresult}
\begin{IEEEproof}
See Appendix~\ref{app:ProofInterpolant}. 
\end{IEEEproof}

Formula~\eqref{eq:interpolation-bound-for-error-in-an-interval} has
the following intuitive explanation: the error measures the deviation
of $\f$ from being straight (linear), and this is
directly related to the convexity/concavity of the function, thus the
presence of the term $\fimax$ to bound the amount of bending.

The other interesting part
of~\eqref{eq:interpolation-bound-for-error-in-an-interval} is the
$h_{i}^{5/2}$ dependence with respect to the interval size or the
local density of knots. Other works in the literature use the weaker
bound that does not require carrying out
integration~\cite{deBoor2001}: $\|\f-\pip\f\|_{\Ltwo(I_{i})}\leq
C\text{\,\ensuremath{\fimax}}h_{i}^{2}$, for some constant $C>0$.

\subsubsection{Best Linear Approximation}

\begin{figure}
\centering\includegraphics[width=84mm]{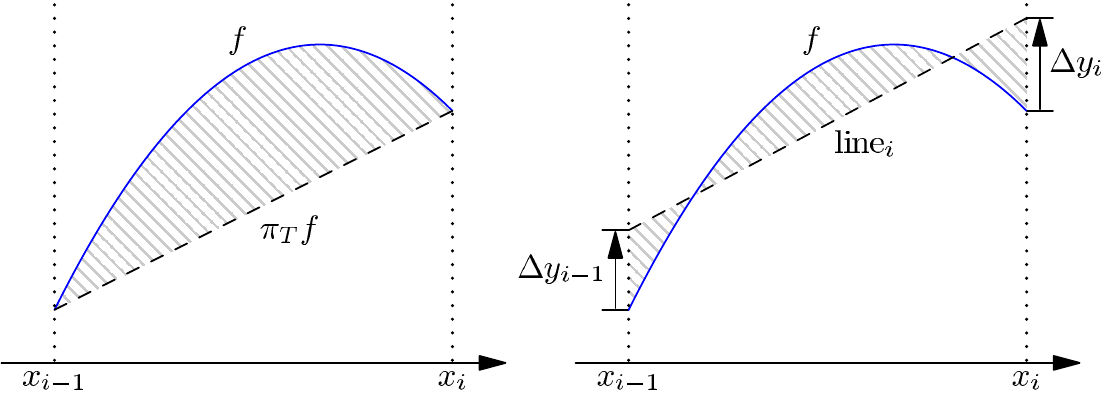}
\caption{\label{fig:errorzone}
Function $\f$ in the subinterval $I_{i}=[x_{i-1},x_{i}]$ and two
linear approximations. On the left, the linear interpolant $\pip\f$
given by~\eqref{eq:LinearInterpolantInInterval}; on the right, a
general linear segment $\mathrm{line}_{i}$ given
by~\eqref{eq:LinearSegment}. $\Delta y_j$ is a signed distance with
respect to $\f(x_j)$.}
\end{figure}

Let us now characterize the error of the orthogonal projection $\ortproj\f$. 
We do so in two steps: 
first we compute the minimum error of a line segment in the interval $I_i$, and then we use an asymptotic analysis to approximate the error of the orthogonal projection.

Stemming from~$\eqref{eq:LinearInterpolantInInterval}$, we can write any
(linear) segment in $I_{i}$ as
\begin{align}
\mathrm{line}_{i}(x) =&\;\bigl(\f(x_{i-1})+\leftDeltaY{i-1}\bigr)\bigl(1-\delta_{i}(x)\bigr)\nonumber\\
&+\bigl(\f(x_{i})+\rightDeltaY
i\bigr)\delta_{i}(x),\label{eq:LinearSegment}
\end{align}
where $\leftDeltaY{i-1}$ and $\rightDeltaY i$ are extra degrees of
freedom (pictured in \figurename~\ref{fig:errorzone}) with respect the
interpolant $\pip\f$ that allow the line segment to better approximate
the function $\f$ in $I_{i}$. 
By computing the optimal values of $\leftDeltaY{i-1}$ and $\rightDeltaY{i}$ we obtain the following result.

\begin{thmresult}\label{thm:MinErrLineSegment} 
The minimum squared $\Ltwo$ error between a given function $f$ and a line segment~\eqref{eq:LinearSegment} in an interval $I_{i}=[x_{i-1},x_{i}]$, of length $h_{i}=x_{i}-x_{i-1}$, 
adopts the expression
\begin{align}
&\min\|\f-\mathrm{line}_{i}\|_{\Ltwo(I_{i})}^{2}=\frac{h_{i}^{5}}{120}(\fsec(\eta))^{2} \nonumber\\
&\quad
-\frac{h_{i}^{5}}{144}\left(\bigl(\fsec(\eta_{0})\bigr)^{2}+\bigl(\fsec(\eta_{1})\bigr)^{2}-\fsec(\eta_{0})\fsec(\eta_{1})\right).\label{eq:strict-best-L2-segment}
\end{align}
for some $\eta,\eta_0,\eta_1$ in $(x_{i-1},x_{i})$.
\end{thmresult}
\begin{IEEEproof}
See Appendix~\ref{app:ProofLineSegment}. 
\end{IEEEproof}
\begin{corollary}\label{thm:MinErrSmallLineSegment}
If $I_{i}$ is sufficiently small so that $\fsec$ is
approximately constant within it, say $\fsec_{I_{i}}$, then
\begin{equation}
\min\|\f-\mathrm{line}_{i}\|_{\Ltwo(I_{i})}^{2} \approx 
\frac{h_{i}^{5}}{720}(\fsec_{I_{i}})^{2}.\label{eq:approximate-best-L2-segment}
\end{equation}
\end{corollary}
\begin{IEEEproof}
Substitute $\fsec(x)\approx\fsec_{I_{i}}$ for $x=\{\eta,\eta_0,\eta_1\}$ in~\eqref{eq:strict-best-L2-segment}.
\end{IEEEproof}

The segments in $\ortproj\f$ may not strictly
satisfy~\eqref{eq:strict-best-L2-segment} because $\ortproj\f$ must
also be continuous across segments. 
However, as the size of the intervals $h_i$ become smaller 
(due to a finer partition~$\part$ of the interval~$I$, i.e., the number of segments $N$ in $\part$ increases) 
we may approximate $\ortproj\f$ in $I_i$ by the best (independently-optimized) segment $\mathrm{line}_{i}$ and, consequently, use Corollary~\ref{thm:MinErrSmallLineSegment} to yield the following result.

\begin{thmresult}\label{thm:MinErrSmallOrtProj} 
The squared $\Ltwo$ error between a given function $f$ and its orthogonal projection $\ortproj\f$ in a \emph{small} interval $I_{i}=[x_{i-1},x_{i}]$, of length $h_{i}=x_{i}-x_{i-1}$, 
is
\begin{equation}
\|\f-\ortproj\f\|_{\Ltwo(I_{i})}^{2} \approx 
\frac{h_{i}^{5}}{720}(\fsec_{I_{i}})^{2}.\label{eq:approximate-best-ortproj}
\end{equation}
where $\fsec_{I_{i}}$ is the approximately constant value that $\fsec$ takes within the interval $I_{i}$.
\end{thmresult}
\begin{IEEEproof}
See Appendix~\ref{app:Convergence}. 
\end{IEEEproof}
In the same asymptotic situation, the linear
interpolant~\eqref{eq:interpolation-bound-for-error-in-an-interval}, (see~\eqref{eq:SquaredErrorLinInterp}) gives a bigger squared error by a
factor of six,
\begin{equation}
\|\f-\pip\f\|_{\Ltwo(I_{i})}^{2} \approx \frac{h_{i}^{5}}{120}(\fsec_{I_{i}})^{2}.\label{eq:approximate-interpolated-segment}
\end{equation}

\subsection{\label{sub:Approximation-to-the-optimal-partition}Approximation
to the Optimal Partition}
Now we give a procedure to compute a suboptimal partition of the target interval $I$ and then derive error estimates for both $\pip \f$ and $\ortproj$ on such a partition and the uniform one.

\subsubsection{Procedure to Obtain a Suboptimal Partition}

Let us consider a partition $\part$ of $I=[a,b]$ with subintervals
$I_{i}=[x_{i-1},x_{i}],\, i=1,\ldots,N$. A suboptimal partition for a
given $N$ is one in which every subinterval has approximately equal
contribution to the total approximation
error~\cite{deBoor2001,Cox2001}, which implies that regions of $\f$
with higher convexity are approximated using more segments than
regions with lower convexity. Let us assume $N$ is large enough so
that $\fsec$ is approximately constant in each subinterval and
therefore the
bound~\eqref{eq:interpolation-bound-for-error-in-an-interval} is
tight. Consequently,
\begin{equation}
\fimax h_{i}^{5/2}\approx C,\label{eq:ErrorEqualization}
\end{equation}
for some constant $C>0$, and the lengths of the subintervals (local
knot spacing~\cite{Cox2001}) should be chosen
$h_{i}\propto\fimax^{-2/5}$, i.e., smaller intervals as $\fimax$
increases. Hence, the local knot distribution or density is
\begin{equation}
lkd(x)\propto|\fsec(x)|^{2/5},\label{eq:limit-optimal-vertex-density}
\end{equation}
so that more knots of the partition are placed in the regions with
larger magnitude of the second derivative.

\begin{figure}
\begin{centering}
\includegraphics[clip,width=84mm]{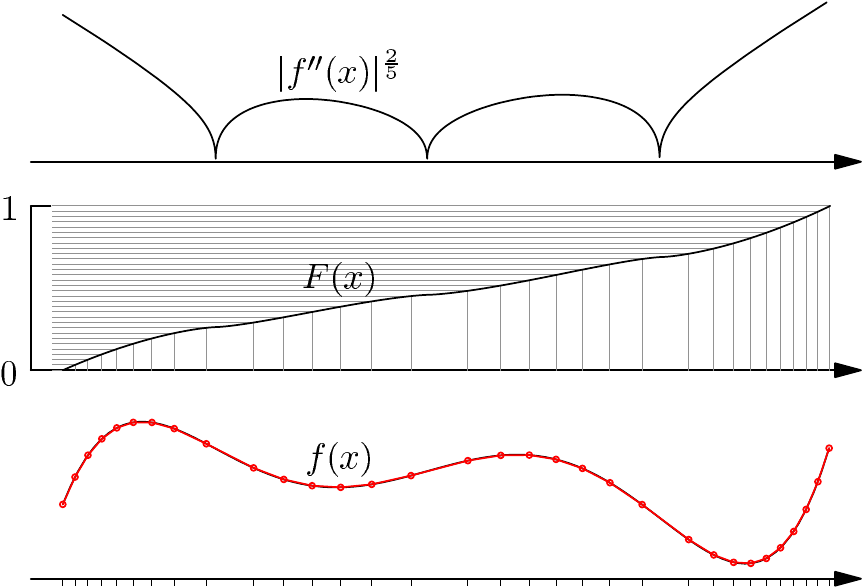}
\par\end{centering}
\caption{\label{fig:optimization}Graphical summary of the proposed 
knot placement technique. 
Top: local knot density~\eqref{eq:limit-optimal-vertex-density} obtained from input function $\f$ (example in \figurename~\ref{fig:pipvsproj});
middle: cumulative knot distribution function $F$ given
by~\eqref{eq:cumulative-knot-distribution} and knots given by the
abscissas corresponding to a uniform partition of the range of $F$, as
expressed by~\eqref{eq:optimized-partition}; 
bottom: approximation of $f$ by CPWL interpolant $\pipOpt\f$ with $N=31$ (32~knots).
In this suboptimal partition, knots are distributed according to the amount of local
convexity/concavity of $\f$ given by the $lkd$ in the middle plot. 
Hence, fewer knots are placed around the zeros of the $lkd$,
which correspond to the less steep regions of $F$. }
\end{figure}

Then, the proposed approximation $\partQuasiOpt$ to the optimal
partition is as follows: $x_{0}=a$, $x_{N}=b$, and take knots
$\{x_{i}\}_{i=1}^{N-1}$ given by
\begin{equation}
F(x_{i})=i/N,\label{eq:optimized-partition}
\end{equation}
where the monotonically increasing function $F:[a,b]\to[0,1]$ is
\begin{equation}
F(x)=\int_{a}^{x}\bigl(\fsec(t)\bigr)^{2/5}\dif t\Bigl/\int_{a}^{b}\bigl(\fsec(t)\bigr)^{2/5}\dif t.\label{eq:cumulative-knot-distribution}
\end{equation}
In this procedure, pictured in \figurename~\ref{fig:optimization}, the range
of $F(x)$ is divided into $N$ contiguous sub-ranges of equal length,
and the values of $x_{i}$ are given by the abscissas corresponding to
the endpoints of the sub-ranges.

\subsubsection{Error Estimates}
With this partition, we can further estimate approximation error bounds in the entire interval based on those of the subintervals.

\begin{thmresult}\label{thm:ErrorSuboptPartPiP}The approximation error of the linear interpolant $\pipQuasiOpt\f$ in the partition $\partQuasiOpt$ given by~\eqref{eq:optimized-partition} in $I=[a,b]$~is
\begin{equation}
\|\f-\pipQuasiOpt\f\|\lesssim\frac{1}{N^{2}\sqrt{120}}\left(\int_{a}^{b}\bigl(\fsec(t)\bigr)^{2/5}\dif t\right)^{5/2}.\label{eq:ApproxL2Boundpip}
\end{equation}
\end{thmresult}

\begin{IEEEproof}
The total squared error for any partition $\part$ is 
the sum of squared errors over all subintervals $I_i$, and by~\eqref{eq:interpolation-bound-for-error-in-an-interval},
\begin{equation}
\|\f-\pip\f\|^{2} \leq \sum_{i=1}^{N}\frac{1}{120}\left(\fimax
h_{i}^{5/2}\right)^{2},\label{eq:SquaredBoundpipL2-prev}
\end{equation}
which, under the condition of error
equalization~\eqref{eq:ErrorEqualization} for the proposed partition
$\partQuasiOpt$, becomes
\begin{equation}
\|\f-\pipQuasiOpt\f\|^{2}\leq\sum_{i=1}^{N}\frac{1}{120}C^{2}=\frac{1}{120}C^{2}N.\label{eq:SquaredBoundpipL2Interval}
\end{equation}
To compute $C$, let us sum $\fimax^{2/5}h_{i}\approx C^{2/5}$ over all
intervals $I_{i}$ and approximate the result using the Riemann
integral:
\begin{equation}
C^{2/5}N \,\stackrel{\eqref{eq:ErrorEqualization}}{\approx}\, 
\sum_{i=1}^{N}\fimax^{2/5}h_{i}\approx\int_{a}^{b}\bigl(\fsec(t)\bigr)^{2/5}\dif t,
\label{eq:CalculationOfConstantC}
\end{equation}
whose right hand side is independent of $N$. 
Substituting~\eqref{eq:CalculationOfConstantC}
in~\eqref{eq:SquaredBoundpipL2Interval} gives the desired approximate bound~\eqref{eq:ApproxL2Boundpip} for
the error in the interval $I=[a,b]$. 
\end{IEEEproof}

Since the approximation error of $\ortproj\f$ in each
interval is roughly proportional to that of $\pip\f$, as shown
in~\eqref{eq:approximate-best-ortproj}
and~\eqref{eq:approximate-interpolated-segment}, the
partition~\eqref{eq:optimized-partition} is also a very good
approximation to the optimal partition for $\ortproj\f$ as the number
of subintervals $N$ increases. 
This is briefly stated next.

\begin{thmresult}\label{thm:ErrorSuboptPartOrtProj}The approximation error of the orthogonal projection $\ortprojQuasiOpt\f$ in the partition $\partQuasiOpt$ given by~\eqref{eq:optimized-partition} in $I=[a,b]$~is
\begin{equation}
\|\f-\ortprojQuasiOpt\f\| \approx \frac{1}{\sqrt{6}}\|\f-\pipQuasiOpt\f\|.\label{eq:ApproxL2BoundOrtProj}
\end{equation}
\end{thmresult}

Both CPWL approximations ($\pip\f$ and $\ortproj\f$) converge to the true function $\f$ at a rate of at least $O(N^{-2})$ (\eqref{eq:ApproxL2Boundpip} and~\eqref{eq:ApproxL2BoundOrtProj}).

We use a similar procedure to derive an estimate error bound for
the uniform partition $\partUnif$ that can be compared to that of the
optimized one. 

\begin{thmresult}\label{thm:ErrorUniformPart}The approximation error of the linear interpolant $\pipUnif\f$ in the uniform partition $\partUnif$ of the interval $I=[a,b]$~is
\begin{equation}
\|\f-\pipUnif\f\| \approx \frac{\left(b-a\right)^{2}}{N^{2}\sqrt{120}}\|\fsec\|.
\label{eq:ApproxL2boundUniformPip}
\end{equation}
\end{thmresult}

\begin{IEEEproof}
For $\partUnif$, we may substitute $h_{i}=(b-a)/N$
in~\eqref{eq:SquaredBoundpipL2-prev} and approximate the result using
the Riemann integral,
\begin{equation*}
\|\f-\pipUnif\f\| \stackrel{\eqref{eq:SquaredBoundpipL2-prev}}{\leq}\frac{h_{i}^{2}}
{\sqrt{120}}\left(\sum_{i=1}^{N}\fimax^{2}h_{i}\right)^{1/2}
\approx\frac{\left(b-a\right)^{2}}{N^{2}\sqrt{120}}\|\fsec\|.
\end{equation*}
\end{IEEEproof}

Thus, we can estimate how much we can expect to benefit from optimizing
a partition by simply dividing~\eqref{eq:ApproxL2Boundpip}
by~\eqref{eq:ApproxL2boundUniformPip}.
Since~\eqref{eq:ApproxL2boundUniformPip} also shows a $O(N^{-2})$ convergence, the profit will not depend on $N$ (assuming $N$ is large enough).

\section{\label{sec:Computational-analysis-and}Computational Analysis}

\begin{algorithm*}[t]
1 \hskip 12pt Determine the subinterval $I_i$ that contains
$x^{\star}$: $i\mid x_{i-1}\leq x^{\star}<x_{i}$

2 \hskip 12pt Find the fractional distance $\delta_i$ to the left
endpoint of $I_i$:
$\delta_i\Leftarrow(x^{\star}-x_{i-1})/(x_{i}-x_{i-1})$

3 \hskip 12pt Return the interpolated value
$\left(1-\delta_i\right)\v(x_{i-1})+\delta_i\,\v(x_{i})$
\caption{(\textbf{CPWL function evaluation}). \label{alg:Given-as-input}Given as input the description of a continuous
piecewise linear function (CPWL) $\v\in\vecspace$ by means of
the the nodal values $\vsamp=\v(x_{0}),\ldots,\v(x_{N})$ at the knots
$\xsamp=x_{0},\ldots,x_{N}$ of the partition $\part$, and an abscissa
$x^{\star}\in\left[x_{0},x_{N}\right)$, return $\v(x^{\star})$.}
\end{algorithm*}

The implementation of the approximations previously discussed is
straightforward (Algorithm~\ref{alg:Given-as-input}) and we only need to distinguish two different cases
depending on whether the partition $\part$ of $I$ has all subintervals
of equal width or not. Although their values are different,
$\pip\f(x)$ and $\ortproj\f(x)$ are qualitatively the same from the
computational viewpoint.

Let us discuss the general case that $\part$ is non-uniform. Line~1 in
Algorithm~\ref{alg:Given-as-input} implies searching in the ordered
array $\xsamp$ for the right index $i$. Since the other steps in
the algorithm do not depend on the size of the input
vectors, Line~1 is the dominant step and makes its run-time
complexity $O(\log N)$. In the particular case that $\part$ is
uniform, no search is needed to determine the index:
$i\Leftarrow\bigl\lfloor\frac{x^{\star}-x_{0}}{x_{N}-x_{0}}N\bigr\rfloor+1$
and the fractional part
$\delta_{i}\Leftarrow\frac{x^{\star}-x_{0}}{x_{N}-x_{0}}N+1-i$. Thus,
the run-time complexity of the uniform case is $O(1)$. There is also
no need to store the full array $\xsamp$ but only its endpoints
$x_{0}$ and $x_{N}$, roughly halving the memory needed for a
non-uniform partition of the same number of intervals.

Consequently, approximations based on a uniform partition are expected
to perform better in usual CPUs or GPUs, computationally-wise, than
those based on non-uniform partitions. However, optimizing the
partition might lead, depending on the specific objective function, to
a reduction in the memory requirements for a desired approximation
error. If memory constraints are relevant or hardware-aided search
routines become available, optimized partitions could become
practical.

\section{\label{sec:Implementation-in-a-GPU}Implementation in a GPU}

The proposed algorithm is simple to implement either in a CPU or in a
GPU. However, there are some implementation details that help further
optimize the performance in a GPU.

Modern GPUs implement a Single Instruction, Multiple Data (SIMD)
execution model in which multiple threads execute exactly the same
instructions (each) over different data. If two threads scheduled
together in the same execution group or \emph{warp} follow different
branches in a conditional statement, both branches are not executed in
parallel but sequentially, thereby halving throughput. In the case of
$\part$ being non-uniform, the (binary) search performed to find the
index $i$ could lead to divergent paths for different threads
evaluating $\v(x)$ for different values of $x$, thus reducing
throughput.

Since the purpose of the proposed algorithm is to save computation at
run-time, it is reasonable to assume that at least $N$, and possibly
$\xsamp$ and $\vsamp$ too, have been determined at compile time.
Then, we can particularize the search code in compile time, unrolling
the loop needed to implement the binary search and, most importantly,
replacing conditional statements with a fixed number of invocations of
the ternary operator, which is implemented in the GPU with a single
instruction and no code divergence~\cite{Pallipuram2012}.

In general, accesses to memory in a GPU need to be arranged so that
neighbouring threads access neighbouring addresses in memory. This
allows the host interface of the GPU to coalesce operations from
several threads into fewer transactions bigger in size. In the
described algorithm, if each thread is evaluating $\v(x)$ for
different values of $x$, it is not guaranteed that the memory access
pattern complies with this restriction. Although some devices provide
caching and isolate the programmer from this issue, not all do.

However, there is a case in a regular computer graphics pipeline in
which a similar problem arises. Texture data need to be accessed many
times in a short period and it is generally impossible to ensure that
neighbouring pixels in screen need to read neighbouring
texels. Consequently, most GPUs do cache texture data and even have
dedicated cache layers for it. We therefore store the arrays of
samples in texture memory rather than in global memory to benefit from
texture caches.

The use of the texture units to access our samples provides another
important benefit. In the usual computer graphics pipeline, GPUs need
to constantly filter texture data. In order to efficiently perform
that task, texture units are equipped with hardware fixed-function
interpolation routines that include linear
interpolation~\cite{Lindholm2008}.  Therefore, we can use them to
speed up line~3 in Algorithm~\ref{alg:Given-as-input}; results in the
next section confirm that the texture filtering unit is indeed faster
than performing interpolation ``manually''.

\section{\label{sec:Results}Experimental Results}

In this section we apply the proposed linearization algorithm to several functions of interest in cybernetics, computer vision and other scientific areas.
The section starts with the analysis of the Gaussian function and demonstrates the proposed technique on an image processing application. 
The section then analyzes two other sample nonlinear functions of interest: the Lorentzian function and the Bessel function of the first kind.

\subsection{Gaussian Function}\label{sec:Gaussian}
The Gaussian function $\f(x)=\exp(-x^{2}/2)/\sqrt{2\pi}$ is
widely used in many applications in every field of science, but it is
of particular interest to us in the context of foreground segmentation 
in video sequences using spatio-temporal nonparametric background
models~\cite{Berjon2013}. To estimate these models, the Gaussian function
needs to be evaluated approximately 1300 times per pixel and input image,
or around 2 billion times per second at modest video quality
(CIF, $352 \times 288$ pixels at 15 fps). Therefore it is crucial to lower
the computing time of the Gaussian function. In the performed
experiments it was enough to approximate the function in the interval
$x\in\left[0,8\right]$ to achieve the required quality. 

\begin{figure}
\begin{centering}
\includegraphics[clip,width=84mm]{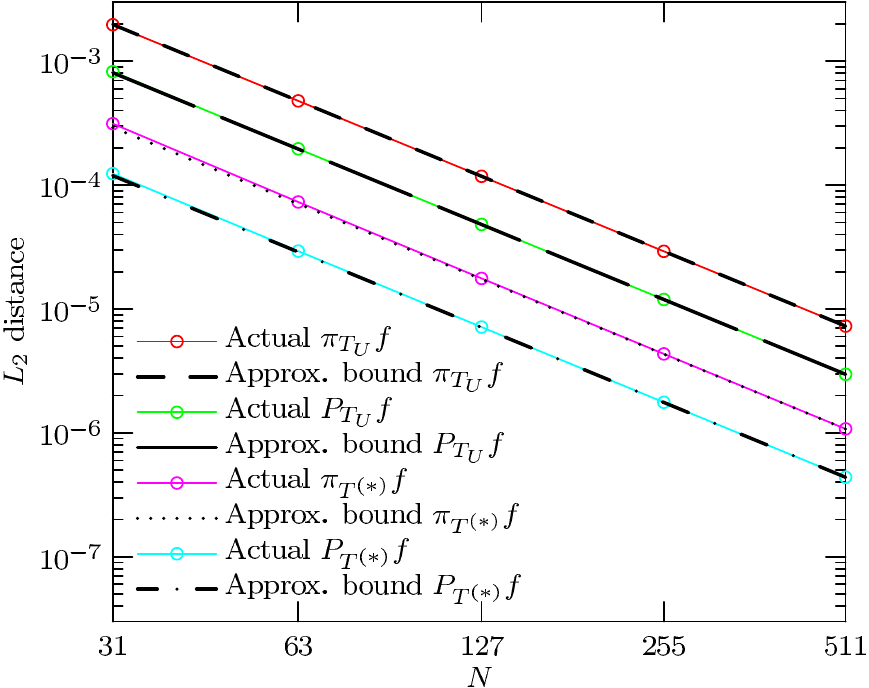}
\par\end{centering}
\caption{\label{fig:l2-distance}$\Ltwo$ distance for different approximations
to the Gaussian function $\f(x)=\exp(-x^{2}/2)/\sqrt{2\pi}$ over the
interval $x\in\left[0,8\right]$.}
\end{figure}

\figurename~\ref{fig:l2-distance} shows $\Ltwo$ distances between the
Gaussian function and the approximations described in
previous sections. It reports actual distances as
well as the approximate and tight upper bounds to the
distances, i.e., Results~\ref{thm:ErrorSuboptPartPiP} to~\ref{thm:ErrorUniformPart}.
It can be observed that the $\Ltwo$ distances between $\pip\f$ and $\f$ using
the uniform and optimized partitions agree well with Results~\ref{thm:ErrorSuboptPartPiP} to~\ref{thm:ErrorUniformPart}.
All curves in \figurename~\ref{fig:l2-distance} have similar slope to that of~\eqref{eq:ApproxL2Boundpip}, i.e., their convergence rate is $O(N^{-2})$, and the ratio
between the distances corresponding to $\ortproj\f$ and $\pip\f$ is
approximately the value $1/\sqrt{6}$ that stems
from Result~\ref{thm:MinErrSmallOrtProj}
and~\eqref{eq:approximate-interpolated-segment}.

\begin{table*}
\caption{\label{tab:Mean-per-evaluation-execution}Mean per-evaluation execution
times (in picoseconds) for the Gaussian function.}
\begin{centering}
\begin{tabular}{|l|l|c|c|c|c|c|}
\hline 
 & {\footnotesize{Number of points ($N+1$)}} & {\footnotesize{32}} &
 {\footnotesize{64}} & {\footnotesize{128}} & {\footnotesize{256}} &
 {\footnotesize{512}}\tabularnewline
\hline 
\hline 
\multirow{3}{*}{{\footnotesize{CPU}}} & {\footnotesize{exact function using }}\texttt{\footnotesize{expf}}{\footnotesize{
from }}\texttt{\footnotesize{<cmath>}}
& \multicolumn{5}{c|}{{\footnotesize{13710}}}\tabularnewline
\cline{2-7} 
 & {\footnotesize{uniform partition}}
 & \multicolumn{5}{c|}{{\footnotesize{1750 - 1780}}}\tabularnewline
\cline{2-7} 
 & {\footnotesize{optimized partition}} & {\footnotesize{8900}} &
 {\footnotesize{11230}} & {\footnotesize{13830}} &
 {\footnotesize{16910}} & {\footnotesize{20120}}\tabularnewline
\hline 
\multirow{6}{*}{{\footnotesize{GPU}}} & {\footnotesize{exact function using }}\texttt{\footnotesize{expf}}{\footnotesize{
from }}\texttt{\footnotesize{<cmath>}}
& \multicolumn{5}{c|}{{\footnotesize{33.3}}}\tabularnewline
\cline{2-7} 
 & {\footnotesize{fast function using
}}\texttt{\footnotesize{\_\_expf}}{\footnotesize{ from the SFU}}
& \multicolumn{5}{c|}{{\footnotesize{14.2}}}\tabularnewline
\cline{2-7} 
 & {\footnotesize{uniform partition, manual interpolation}}
 & \multicolumn{5}{c|}{{\footnotesize{19.7 - 19.8}}}\tabularnewline
\cline{2-7} 
 & {\footnotesize{uniform partition, hardware interpolation}}
 & \multicolumn{5}{c|}{{\footnotesize{7.8 - 7.9}}}\tabularnewline
\cline{2-7} 
 & {\footnotesize{optimized partition, manual interpolation}} &
 {\footnotesize{142.0}} & {\footnotesize{161.9}} &
 {\footnotesize{181.8}} & {\footnotesize{203.1}} &
 {\footnotesize{224.1}}\tabularnewline
\cline{2-7} 
 & {\footnotesize{optimized partition, hardware interpolation}} &
 {\footnotesize{122.7}} & {\footnotesize{142.9}} &
 {\footnotesize{163.2}} & {\footnotesize{188.3}} &
 {\footnotesize{211.1}}\tabularnewline
\hline 
\end{tabular}
\par\end{centering}
\end{table*}

Table~\ref{tab:Mean-per-evaluation-execution} shows mean processing
times per evaluation both in the CPU (sequentially, one core only) and
in the GPU. All execution times have been measured in a computer 
equipped with an Intel Core~i7-2600K processor, 16~GiB
RAM and an NVIDIA GTX~580 GPU. We have exercised the utmost care to
ensure that all different strategies achieve the same utilization of
the GPU so that measurements are fair.

\begin{figure}
\begin{centering}
\includegraphics[clip,width=84mm]{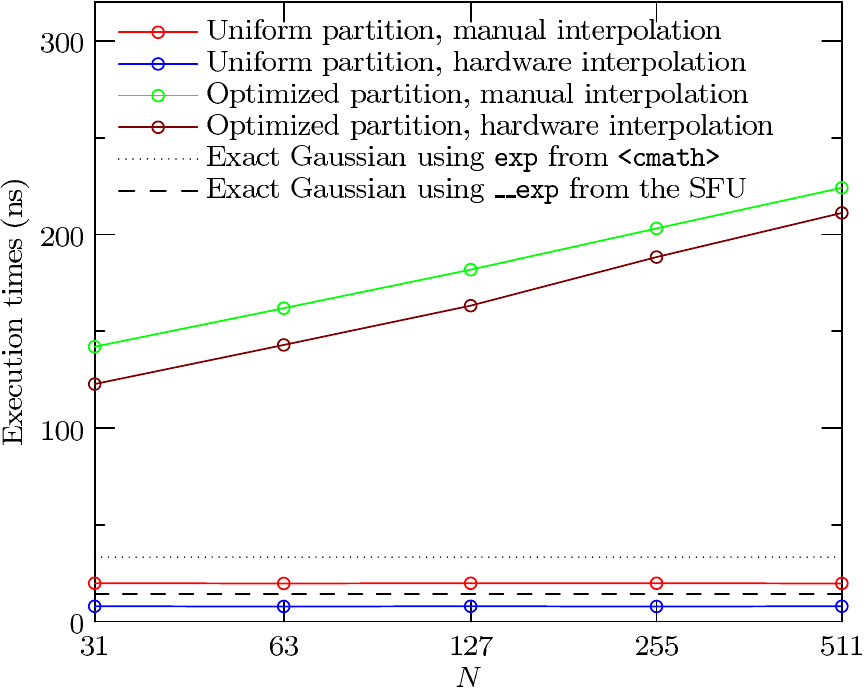}
\par\end{centering}
\caption{\label{fig:times}Mean per-evaluation execution times of all
the proposed variants on a GPU (in picoseconds). The graph clearly shows 
the $O(1)$ and $O(\log N)$ dependencies on the number of intervals $N$ 
for uniform and non-uniform partitions, respectively.}
\end{figure}

We compare the proposed algorithm against the exact Gaussian
implemented using the exponential function of the C standard library (CPU or GPU)
and against a fast Gaussian implemented using the hardware-accelerated
low-precision exponential function from the SFU of the
GPU~\cite{NVIDIACorporation2012}.  
In either hardware (CPU or GPU), there are no separate
measurements for $\pip\f(x)$ and $\ortproj\f(x)$ because they only
differ in value but not in their implementation. 
As was expected from the computational analysis, execution times for the uniform partition
are constant, whereas they are heavily dependent on~$N$ for the optimized partition
(see \figurename~\ref{fig:times}).

In the CPU case, the optimized partition approach is faster than the C standard library implementation if $N<128$, but the speed-up gain is smaller than 
that of the uniform partition approach, which is 8 times faster than 
the C standard library (regardless of $N$).

In the GPU case, we have also measured the proposed
algorithm both with and without the texture filtering unit 
(last four rows of Table~\ref{tab:Mean-per-evaluation-execution}) to prove
that the former option is indeed faster than coding the interpolation
explicitly, thanks to its dedicated circuitry.
The proposed strategy, using a uniform partition (7.9 ps), solidly outperforms
both the exact (33.3 ps) and fast SFU (14.2 ps) Gaussians, 
by approximate speed-up factors of $\times 4$ and $\times 2$, respectively, 
even though the latter is hardware-accelerated.  
However, we must stress that this is a synthetic benchmark; 
in a real application such as the one described
in the next section, there may exist limitations in memory bandwidth to 
access operands. Moreover, many resources such as shared memory, cache and 
hardware registers are being used for other tasks too. This can (and does) 
affect the execution time of each of the different implementation strategies.

\subsubsection*{\label{sec:segmentation}Foreground Segmentation (Sample Application)}

\begin{figure}
\begin{minipage}[t]{0.5\columnwidth}%
\begin{center}
\includegraphics[width=0.996\columnwidth,height=0.747\columnwidth]{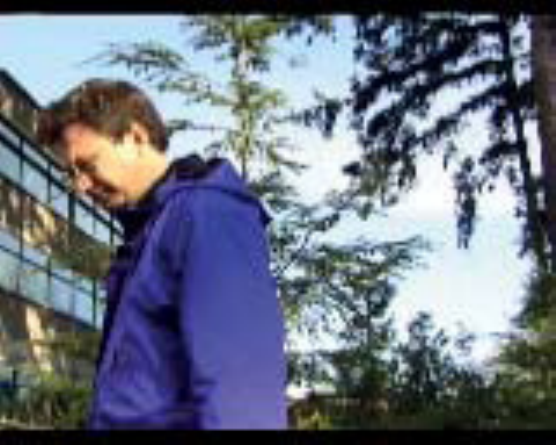}\\
{\footnotesize{}(a)}
\par\end{center}{\footnotesize \par}

\begin{center}
\includegraphics[width=0.996\columnwidth,height=0.747\columnwidth]{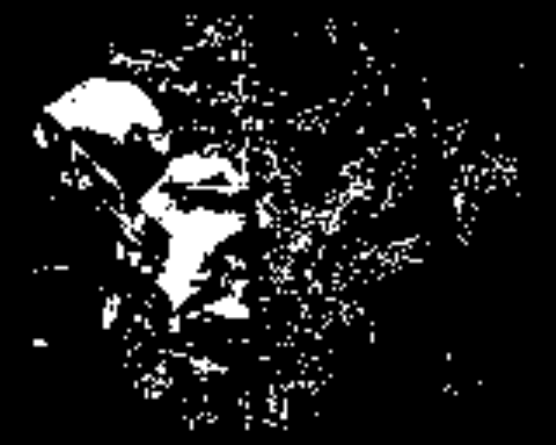}\\
{\footnotesize{}(c)}
\par\end{center}{\footnotesize \par}

\begin{center}
\includegraphics[width=0.996\columnwidth,height=0.747\columnwidth]{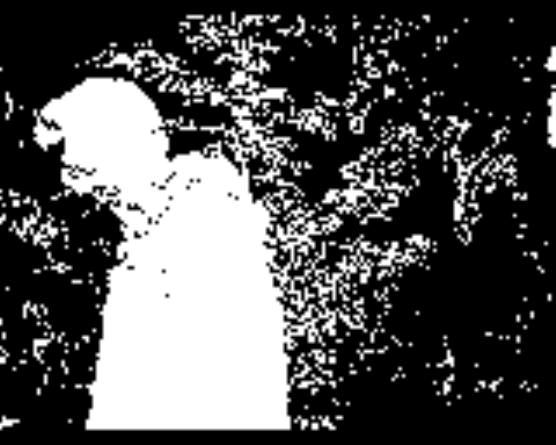}\\
{\footnotesize{}(e)}
\par\end{center}%
\end{minipage}%
\begin{minipage}[t]{0.5\columnwidth}%
\begin{center}
\includegraphics[width=0.996\columnwidth,height=0.747\columnwidth]{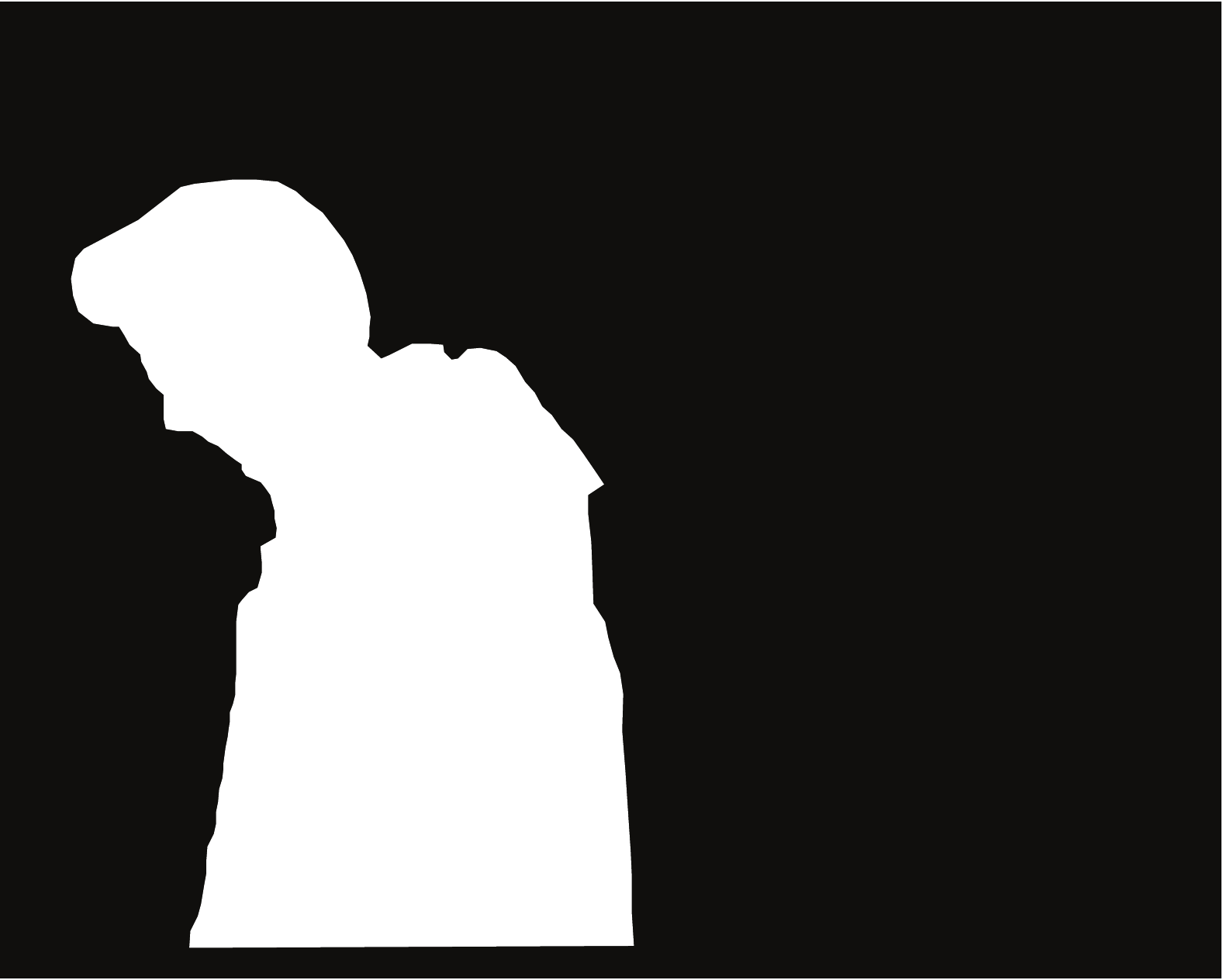}\\
{\footnotesize{}(b)}
\par\end{center}{\footnotesize \par}

\begin{center}
\includegraphics[width=0.996\columnwidth,height=0.747\columnwidth]{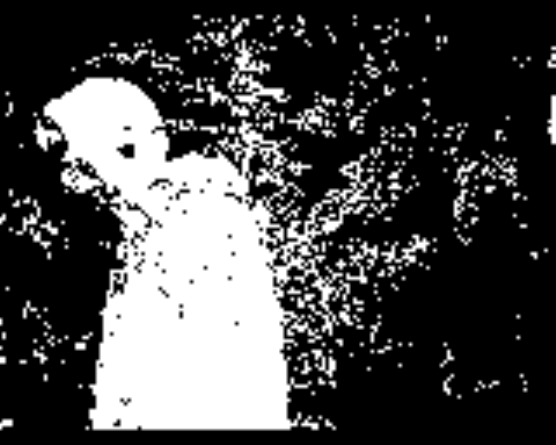}\\
{\footnotesize{}(d)}
\par\end{center}{\footnotesize \par}

\begin{center}
\includegraphics[width=0.996\columnwidth,height=0.747\columnwidth]{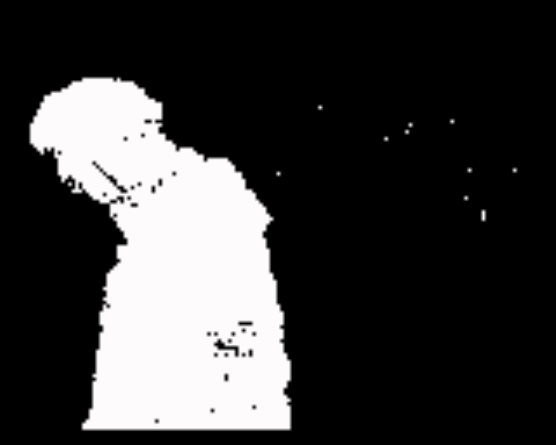}\\
{\footnotesize{}(f)}
\par\end{center}%
\end{minipage}

\caption{\label{fig:mod-methods}Results from several different foreground segmentation methods for a complex dynamic scene. (a) Original scene; (b) Ground truth; (c) Temporal Median Filter; (d) Running Gaussian Average; (e) Mixture of Gaussians; (f) Nonparametric modeling.}
\end{figure}

To obtain high-quality detections in complex situations (e.g.,
dynamic background, illumination changes, etc.), multimodal
moving object detection strategies are commonly used because
they are able to model the multiple pixel states in such
situations. Among multimodal strategies, nonparametric methods
have shown to provide the best quality
detections~\cite{Cuevas2013-IVC} because, in contrast to other
algorithms, they do not assume that the pixel values conform to a
particular distribution, but obtain instead probabilistic models
from sets of recent samples. 

\figurename~\ref{fig:mod-methods} illustrates 
this by comparing the results from different strategies against a 
difficult scene featuring a dynamic background.
The results of two unimodal background modeling methods 
are shown in the second row of this figure: 
the Temporal Median Filter (e.g., \cite{cucchiara2003detecting}) and 
the Running Gaussian Average (e.g., \cite{tang2007background});
whereas the segmentations resulting from 
the application of two multimodal modeling approaches
are depicted in the third row: 
an improved version of the Mixture of Gaussians method~\cite{cuevas2008new} and 
the spatio-temporal nonparametric-based background modeling in~\cite{Berjon2013}.

\begin{figure}
\begin{minipage}[t]{0.5\columnwidth}%
\begin{center}
\includegraphics[width=0.996\columnwidth,height=0.747\columnwidth]{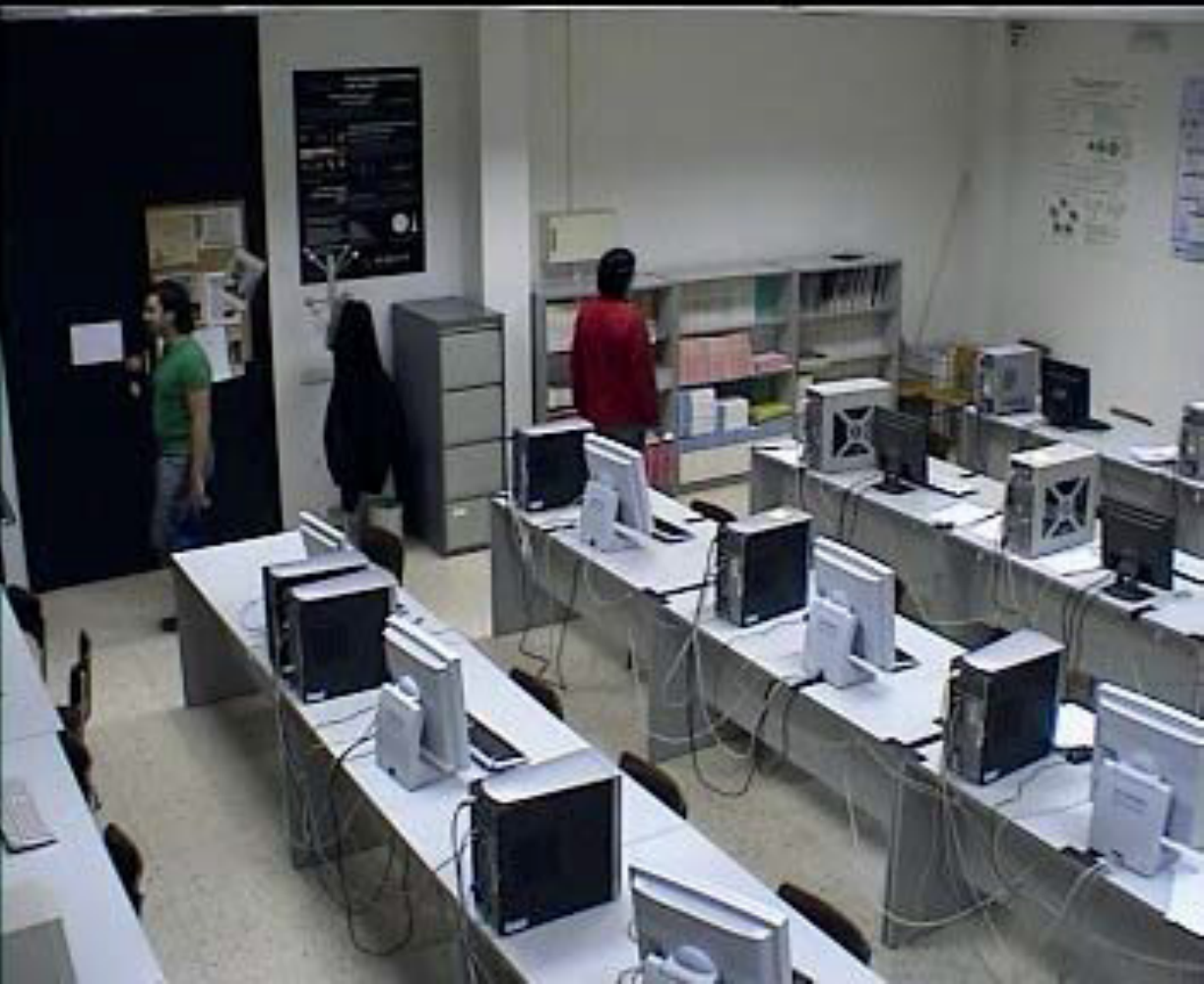}\\
{\footnotesize{}(a)}
\par\end{center}{\footnotesize \par}

\begin{center}
\includegraphics[width=0.996\columnwidth,height=0.747\columnwidth]{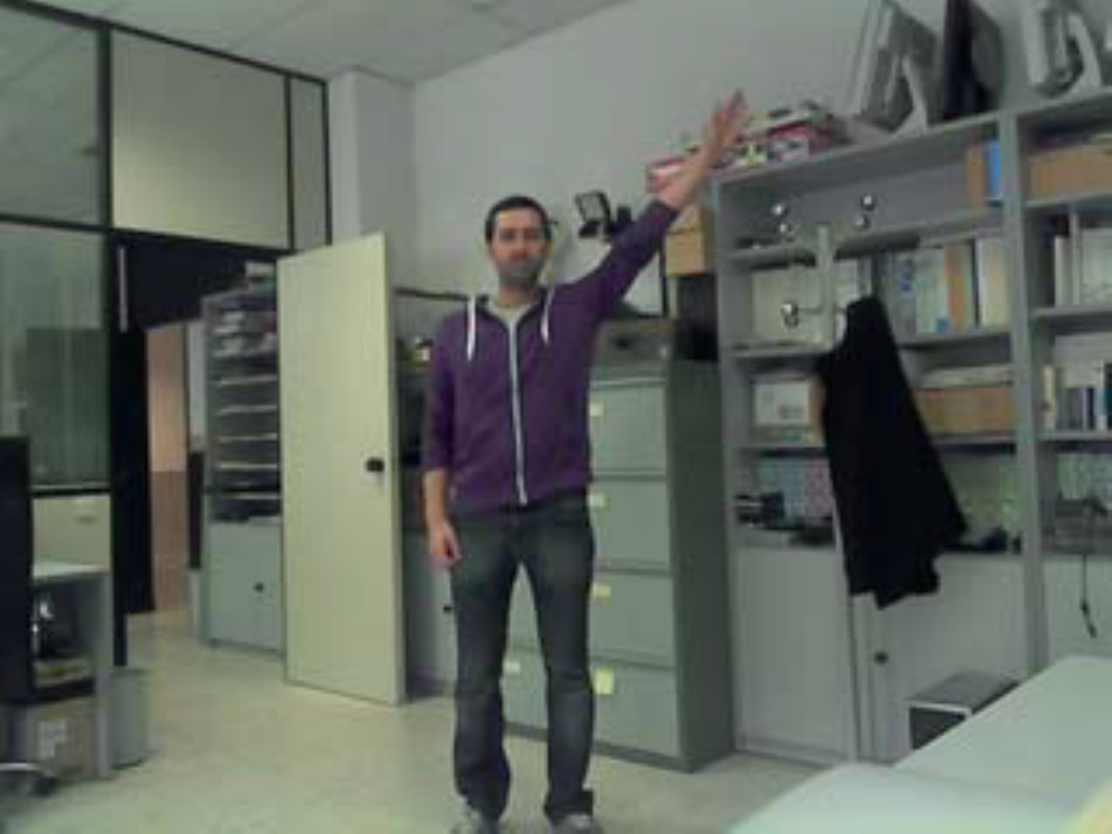}\\
{\footnotesize{}(c)}
\par\end{center}%
\end{minipage}%
\begin{minipage}[t]{0.5\columnwidth}%
\begin{center}
\includegraphics[width=0.996\columnwidth,height=0.747\columnwidth]{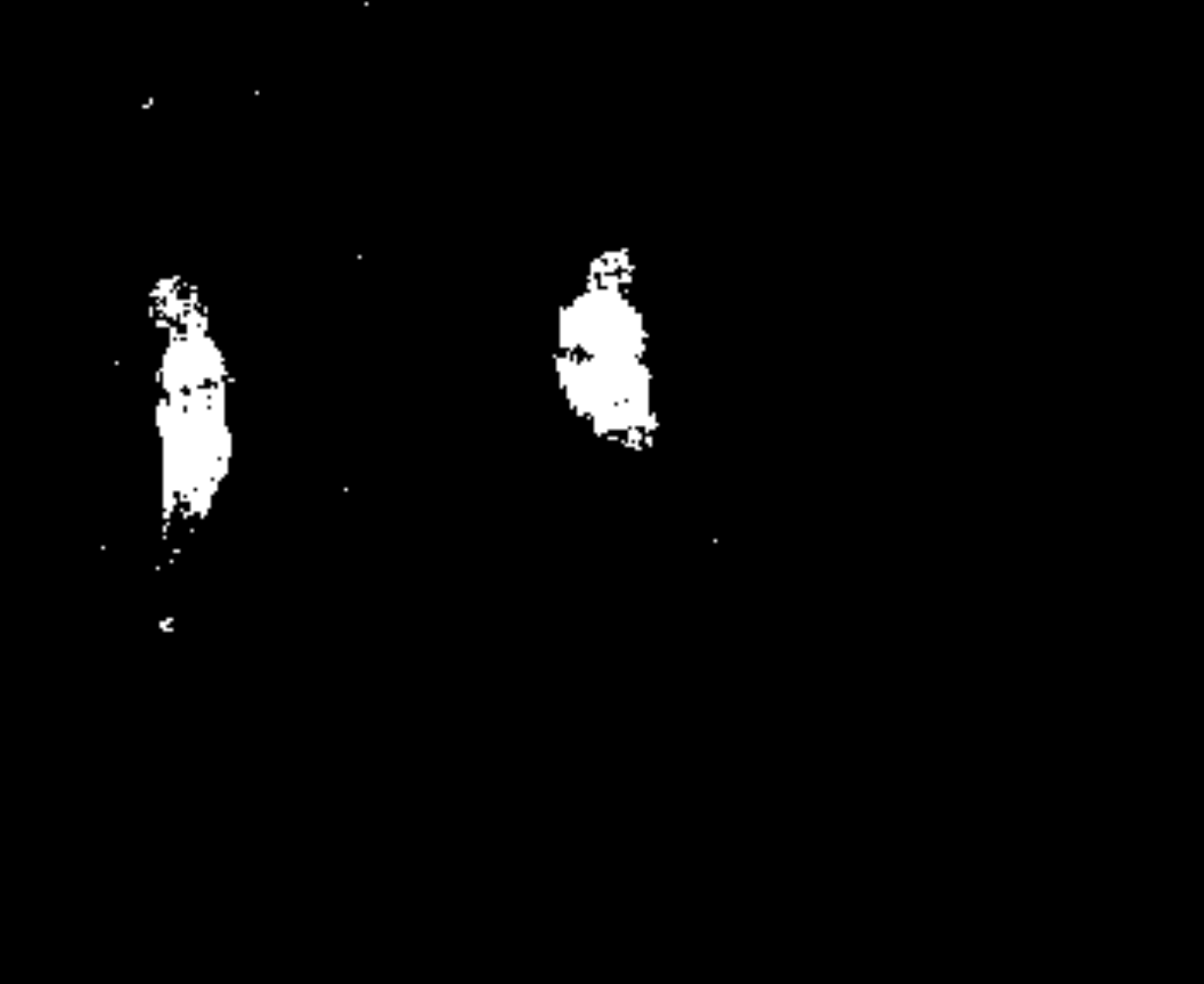}\\
{\footnotesize{}(b)}
\par\end{center}{\footnotesize \par}

\begin{center}
\includegraphics[width=0.996\columnwidth,height=0.747\columnwidth]{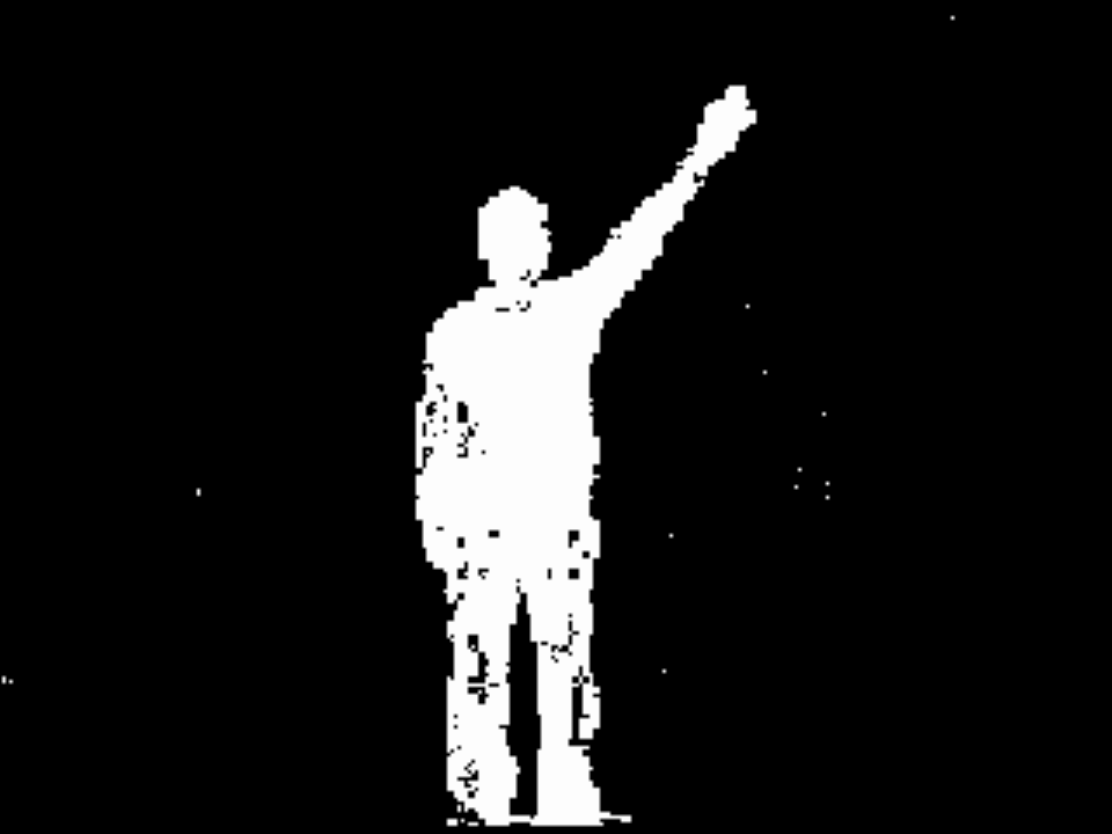}\\
{\footnotesize{}(d)}
\par\end{center}%
\end{minipage}

\caption{\label{fig:modresults}Sample detections using nonparametric
modeling with CPWL-approximated Gaussian kernels. 
Subfigures (a) and (b) show the original scene and results of the segmentation, respectively, for the tracking application; subfigures (c) and (d) show the original scene and results of the segmentation, respectively, for the video-based interface application.}
\end{figure}

To improve the quality of the results in sequences containing non-static 
background regions, and/or captured with portable cameras, 
most recent nonparametric methods use spatio-temporal 
reference data~\cite{sheikh2005bayesian}.
However, these strategies involve huge memory and computational
costs, since they need to evaluate millions of Gaussian kernels per
image. Therefore, lowering the cost of evaluating the 
Gaussian function has a direct impact in the performance of the 
algorithm.
To demonstrate this, we have implemented our proposed technique in a 
nonparametric background modeling method~\cite{Berjon2013}. 
Table~\ref{tab:modtimes} shows mean per-pixel 
(in order to make measurements resolution-independent) processing times 
for the three best performing options in Section~\ref{sec:Results}.
Processing times are not proportional to those in 
Table~\ref{tab:Mean-per-evaluation-execution} because the model needs to do
many other things aside from evaluating Gaussians, and access to reference
data is bandwidth-limited. 
However, our proposed technique still outperforms any of the alternatives:
the exact Gaussian by 27\% and the fast SFU Gaussian by 12\%,
while the final segmentation results of the method are the same in all three cases, 
as depicted in \figurename~\ref{fig:modresults}.

\begin{table}
\caption{\label{tab:modtimes}Mean per-pixel processing times (in nanoseconds) for the nonparametric background modeling technique.}
\begin{center}
\begin{tabular}{lc}
\hline 
{\footnotesize{}Uniform partition, hardware interpolation} & {\footnotesize{}286}\tabularnewline
{\footnotesize{}Fast Gaussian using }\texttt{\footnotesize{}\_\_expf}{\footnotesize{}
from the SFU} & {\footnotesize{}324}\tabularnewline
{\footnotesize{}Exact Gaussian using }\texttt{\footnotesize{}expf}{\footnotesize{}
from }\texttt{\footnotesize{}<cmath>} & {\footnotesize{}391}\tabularnewline
\hline 
\end{tabular}
\end{center}
\end{table}

\subsection{Lorentzian Function}

The performance of the approximation technique has also been tested on other functions.
For example, the Lorentzian function with peak at $x_0$ and width $\gamma$ is
\begin{equation}\label{eq:LorentzFunction}
\mathcal{L}(x; x_0,\gamma) = \frac1{\pi}  \frac{\gamma}{(x - x_0)^2 + \gamma^2}
\end{equation}
As a probability distribution, \eqref{eq:LorentzFunction} is known as the Cauchy distribution.
It is an important function in Physics since it solves the differential equation describing forced resonance. 
In such case, $x_0$ is the resonance frequency and $\gamma$ depends on the damping of the oscillator (and is inversely proportional to the $Q$ factor, a measure of the sharpness of the resonance).

\begin{figure}
\begin{centering}
\includegraphics[clip,width=84mm]
{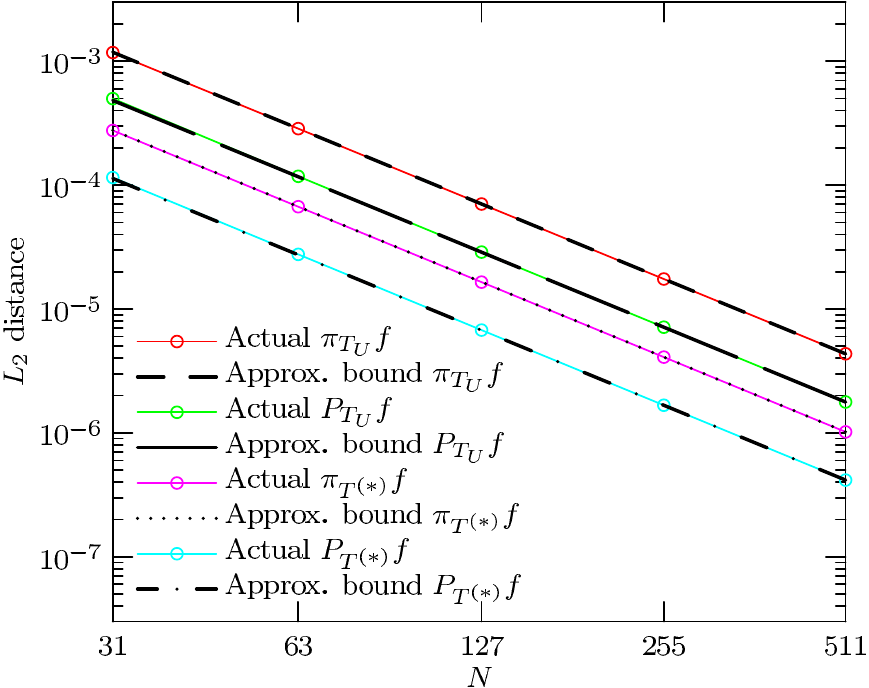}
\par\end{centering}
\caption{\label{fig:GraphicalErrorsCauchy}$\Ltwo$ distance for different CPWL
approximations to the standard Cauchy distribution $f(x) = 1/((1+x^2)\pi)$ in the interval $x\in[0,6]$.}
\end{figure}
\figurename~\ref{fig:GraphicalErrorsCauchy} reports
the $\Ltwo$ distances between $\mathcal{L}(x;0,1)$~\eqref{eq:LorentzFunction} and the CPWL approximations described in previous sections, in the interval $x\in[0,20]$.
The measured errors agree well with the predicted approximate error bounds, 
showing the expected $O(N^{-2})$ trend.

We have measured the mean per-evaluation execution times of the exact
function both in the CPU (576~ps) and in the GPU, both using regular
and reduced-precision accelerated division (19.5~ps and 9.2~ps,
respectively). The evaluation times of our proposed strategy coincide
with those of the Gaussian function
(Table~\ref{tab:Mean-per-evaluation-execution}) because the processing
time of the CPWL approximation does not depend on the function values.

As expected from the lower complexity of this function (compared to that of Sec.~\ref{sec:Gaussian}),
which only involves elementary arithmetic operations, the advantage of our
proposal vanishes in the CPU because the operations needed to manually
perform interpolation, together with the two values that need to be
fetched from memory, are actually more than what the direct evaluation
requires.  However, note that in the GPU our proposal still remains
competitive due to these operations being carried out by the dedicated
circuitry of the texture unit.

\subsection{Bessel Function}

\begin{figure}
\centering{}
\includegraphics[clip,width=84mm]
{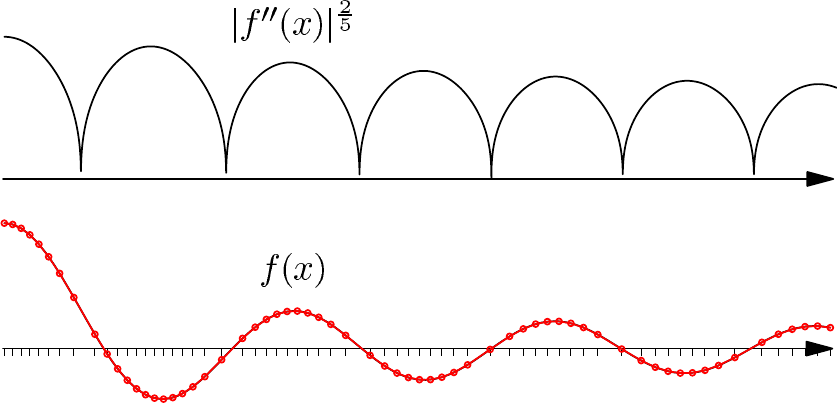}
\caption{\textbf{\label{fig:BesselJ0}}Suboptimal partition for the Bessel function $f(x)=J_0(x)$ in the interval $x\in[0,20]$. 
Top: local knot density ($lkd$) corresponding to $f$; 
bottom: CPWL interpolant $\pipOpt f$ with $N=63$ ($64$ knots) overlaid on function $f$.}
\end{figure}

\begin{figure}
\begin{centering}
\includegraphics[clip,width=84mm]
{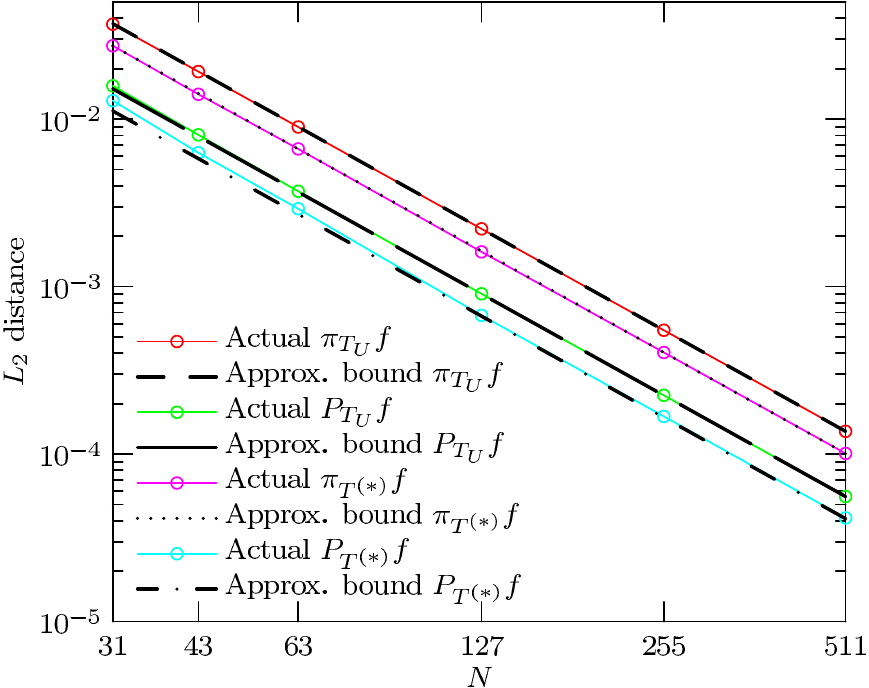}
\par\end{centering}
\caption{\label{fig:GraphicalErrorsBesselJ0}$\Ltwo$ distance for different CPWL
approximations to the Bessel function of the first kind $f(x)=J_0(x)$ in the interval $x\in[0,20]$.}
\end{figure}

Approximating the Bessel function of the first kind $J_0(x)$ is more challenging than approximating the Gaussian or Lorentzian (bell-shaped) functions because it combines both nearly flat and oscillatory parts (see \figurename~\ref{fig:BesselJ0}).
\figurename~\ref{fig:GraphicalErrorsBesselJ0} presents the $\Ltwo$ approximation errors for $J_0(x)$ using CPWL functions in the interval $x\in[0,20]$.
For $N\geq63$ the measured errors agree well with the predicted approximate
error bounds, whereas for $N<63$ the measured errors slightly differ from
the predicted ones (specifically in the optimized partition) because
in these cases the scarce number of linear segments does not properly
represent the oscillations.

A sample optimized partition and the corresponding CPWL interpolant
$\pipOpt f$ is also represented in \figurename~\ref{fig:BesselJ0}. The knots
of the partition are distributed according to the local knot density (see \figurename~\ref{fig:BesselJ0}, Top), 
accumulating in the regions of high oscillations (excluding the places around
the zeros of the $lkd$).

This function is more complex to evaluate in general because it does
not have a closed form. However, our approximation algorithm works
equally well on it.  We have measured the mean per-evaluation
execution times in the CPU (39~ns using the POSIX extensions of the
GNU~C Library and 130~ns---only double-precision provided---using the 
GNU~Scientific Library) and in the GPU (78~ps using the CUDA~\cite{NVIDIACorporation2012} standard library). 
We have also measured the execution time of the first term in the asymptotic expansion~\cite[Eq.~9.57a]{BronshteinHandbook}, valid for $x\gg 1/4$,
$
J_0(x) = \sqrt{\frac{2}{\pi x}}\left[\cos \left(x-\frac{\pi}{4}\right) + O\left(\frac1{x}\right)\right].
$
In the GPU the evaluation takes 42~ps using reduced-precision accelerated functions from
the SFU for the cosine and multiplicative inverse of the square root. 
Despite this approximation having a non-customizable error (decreasing as $x$ increases) and
using the fastest available implementation (SFU) of the operations
involved, our strategy still outperforms it by a sizable margin. This
clearly illustrates that the more complex the function to be
evaluated is, the greater the advantage of our proposal.


\section{\label{sec:Conclusions}Conclusions}

We have developed a fast method to numerically evaluate any continuous
mathematical function in a given interval using
simpler continuous piecewise linear (CPWL) functions and
the built-in routines present in the texture units of modern GPUs. Our
technique allows real-time implementation of demanding computer vision 
and cybernetics applications that use such mathematical functions.

For this purpose, we analyzed the CPWL approximations given by
the linear interpolant and the $\Ltwo$ orthogonal projection of a
function.  We carried out a detailed error analysis in the $\Ltwo$
distance to seek a nearly optimal design of both approximations. In
the practical asymptotic case of a large number of
subintervals $N$, we used error equalization to achieve a
suboptimal design (partition $\partOpt$) and derived a tight bound on the
approximation error for the linear interpolant, showing a $O(N^{-2})$
convergence rate that was confirmed by experimental results. The
$\Ltwo$ orthogonal projection can only but improve upon the results of
the linear interpolant, resulting in a gain factor of
$\sqrt{6}$.

We discussed the computational complexity and the implementation in a
GPU of the numerical evaluation of both CPWL approximations.  Our
experimental results show that our technique can outperform both the
quality and the computational cost of previous similar approaches. In
particular, the fastest strategy consists of using the texture units
in the GPU to evaluate either of the CPWL approximations defined over
a uniform partition. This is normally faster than performing exact
function evaluations, even when using the reduced-precision accelerated
implementation of the SFUs in a GPU, or evaluating the proposed linear
approximations without the assistance of the texture units.
In practice, the number of subintervals $N$ to be considered for
representing the nonlinear function can be decided on its own or based
on other considerations such as a target approximation error, speed or
memory constraints.

Although mathematically sound, the strategies based on a suboptimal
(non-uniform) partition are not practical 
to implement in current CPU/GPU architectures
due to the high cost incurred to find the
subinterval that contains the given point of evaluation. 
Nevertheless, this opens future research paths to explore 
practical implementations of such approaches using specialized hardware.

\appendices

\section{Proof of Result~\ref{thm:LinInterpolantErrBound} (Linear Interpolant Error)}
\label{app:ProofInterpolant}
Recall one of the theorems on
interpolation errors~\cite{Cheney2012}. Let $\f$ be a function in
$C^{n+1}[a,b]$, and let $p$ be a polynomial of degree $n$ or less that
interpolates the function $\f$ at $n+1$ distinct points
$x_{0},x_{1},\ldots,x_{n}\in[a,b]$.  Then, for each $x\in[a,b]$ there
exists a point $\xi_{x}\in[a,b]$ for which
\begin{equation}
\f(x)-p(x)=\frac{1}{(n+1)!}\f^{(n+1)}(\xi_{x})\prod_{k=0}^{n}(x-x_{k}).\label{eq:ThmErrorPolInterp}
\end{equation}

In the interval $I_{i}$, the linear interpolant $\pip\f$ is given by~\eqref{eq:LinearInterpolantInInterval}.
Since $\pip\f$ interpolates
the function $\f$ at the endpoints of $I_{i}$, we can apply
theorem~\eqref{eq:ThmErrorPolInterp} (with $n=1$); hence, the
approximation error solely depends on $\fsec$ and $x$, but not on $\f$
or~$\f^{\prime}$:
\begin{equation}
\f(x)-\pip\f(x)=-\frac{1}{2}\fsec(\xi_{x})(x-x_{i-1})(x_{i}-x).\label{eq:interpolant-difference}
\end{equation}

Let us integrate the square of~\eqref{eq:interpolant-difference} over
the interval $I_{i}$,
\[
\|\f-\pip\f\|_{\Ltwo(I_{i})}^{2} = \int_{I_{i}}\frac{\bigl(\fsec(\xi_{x})\bigr)^{2}}{4}\bigl((x-x_{i-1})(x_{i}-x)\bigr)^{2}\,\dif x.
\]

Next, to simplify the previous integral, let us use the first mean
value theorem for integration, which states that if
$u:[A,B]\rightarrow\R$ is a continuous function and $v$ is an
integrable function that does not change sign on the interval $(A,B)$,
then there exists a number $\eta\in(A,B)$ such that
\begin{equation}
\int_{A}^{B}u(x)v(x)\,\dif x=u(\eta)\int_{A}^{B}v(x)\,\dif x.\label{eq:FirstMVTIntegration}
\end{equation}

Since $(x-x_{i-1})\geq0$ and $(x_{i}-x)\geq0$ for all $x\in I_{i}$, let us
apply~\eqref{eq:FirstMVTIntegration} to compute
\begin{align}
\|\f-\pip\f\|_{\Ltwo(I_{i})}^{2} &=\frac{1}{4}\bigl(\fsec(\eta)\bigr)^{2}\int_{I_{i}}\bigl((x-x_{i-1})(x_{i}-x)\bigr)^{2}\,\dif x \nonumber\\
&=\frac{h_{i}^{5}}{120}\bigl(\fsec(\eta)\bigr)^{2},\label{eq:SquaredErrorLinInterp}
\end{align}
for some $\eta\in(x_{i-1},x_{i})$. Finally, if
$\fimax=\max_{\eta\in I_{i}}|\fsec(\eta)|$, it is
straightforward to derive the $\Ltwo$ error bound~\eqref{eq:interpolation-bound-for-error-in-an-interval} from the square root of~\eqref{eq:SquaredErrorLinInterp}.

\section{Proof of Result~\ref{thm:MinErrLineSegment} (Line Segment Minimum Error)}
\label{app:ProofLineSegment}
The approximation error corresponding to the line segment~\eqref{eq:LinearSegment} is
\begin{equation}
\label{eq:f_minus_line}
\f(x)-\mathrm{line}_{i}(x)=\f(x)-\pip\f(x)-(\pip\Delta y)(x)
\end{equation}
where, by analogy with the
form~$\eqref{eq:LinearInterpolantInInterval}$ of $\pip\f,$ we defined
$(\pip\Delta
y)(x)=\leftDeltaY{i-1}\bigl(1-\delta_{i}(x)\bigr)+\rightDeltaY
i\delta_{i}(x)$.

The proof proceeds by computing the optimal values of $\leftDeltaY{i-1}$ and $\rightDeltaY{i}$ that minimize the squared $\Ltwo$ error over the interval $I_{i}$:
\begin{align}
\|\f-\mathrm{line}_{i}\|_{\Ltwo(I_{i})}^{2} & \stackrel{\eqref{eq:f_minus_line}}{=} \|\f-\pip\f\|_{\Ltwo(I_{i})}^{2}+\|\pip\Delta
 y\|_{\Ltwo(I_{i})}^{2}\nonumber\\ & \quad
 -2\inner{\f-\pip\f}{\pip\Delta
 y}_{\Ltwo(I_{i})}.\label{eq:ExpandedSquaredNorm}
\end{align}
The first term is given in~\eqref{eq:SquaredErrorLinInterp}. 
The second term is
\[
\|\pip\Delta y\|_{\Ltwo(I_{i})}^{2} = \frac{h_{i}}{3}\left((\leftDeltaY{i-1})^{2}+(\rightDeltaY
 i)^{2}+\leftDeltaY{i-1}\rightDeltaY i\right),
\]
and the third term is, applying~\eqref{eq:FirstMVTIntegration} 
and the change of variables $t=\delta_{i}(x)$ to evaluate the resulting integrals,
\begin{align*}
&-2\inner{\f-\pip\f}{\pip\Delta y}_{\Ltwo(I_{i})}\\
 &
 =\leftDeltaY{i-1}\int_{I_{i}}\fsec(\xi_{x})(x-x_{i-1})(x_{i}-x)\bigl(1-\delta_{i}(x)\bigr)\,\dif
 x\\ & \quad +\rightDeltaY
 i\int_{I_{i}}\fsec(\xi_{x})(x-x_{i-1})(x_{i}-x)\delta_{i}(x)\,\dif x\\
& =\frac{h_{i}^{3}}{12}\bigl(\leftDeltaY{i-1}\fsec(\eta_{0})+\rightDeltaY
i\fsec(\eta_{1})\bigr),
\end{align*}
for some $\eta_{0}$ and $\eta_{1}$ in $(x_{i-1},x_{i})$.

Substituting previous results in~\eqref{eq:ExpandedSquaredNorm}, 
\begin{align*}
\|\f-\mathrm{line}_{i}\|_{\Ltwo(I_{i})}^{2} &= 
\frac{h_{i}}{3}\Bigl( (\leftDeltaY{i-1})^{2}+(\rightDeltaY i)^{2}+\leftDeltaY{i-1}\rightDeltaY i \Bigr)\\
&\quad +\frac{h_{i}^{3}}{12}\bigl(\leftDeltaY{i-1}\fsec(\eta_{0})+\rightDeltaY i\fsec(\eta_{1})\bigr)\\
&\quad + \frac{h_{i}^{5}}{120}\bigl(\fsec(\eta)\bigr)^{2}.
\end{align*}
We may now find the line segment that minimizes the distance to $\f$
by taking partial derivatives with respect to $\leftDeltaY{i-1}$ and
$\rightDeltaY i$, setting them to zero and solving the corresponding
system of equations. Indeed, the previous error is quadratic in
$\leftDeltaY{i-1}$ and $\rightDeltaY i$, and attains its minimum at
\begin{equation}\label{eq:DeltaYMinimizers}
\begin{cases}
\leftDeltaY{i-1} & = \bigl(\fsec(\eta_{1})-2\fsec(\eta_{0})\bigr)h_{i}^{2}/12,\\
\quad\rightDeltaY i & = \bigl(\fsec(\eta_{0})-2\fsec(\eta_{1})\bigr)h_{i}^{2}/12.
\end{cases}
\end{equation}
The resulting minimum squared distance is
\eqref{eq:strict-best-L2-segment}.

\section{Proof of Result~\ref{thm:MinErrSmallOrtProj} (Asymptotic Analysis)}
\label{app:Convergence}
By the triangle inequality, 
the jump discontinuity at $x=x_{i}$ between two adjacent independently-optimized segments is
\begin{align}
\left|\rightDeltaY i^{-}-\leftDeltaY i^{+}\right| 
\nonumber & \leq\frac{h_{i}^{2}}{12}\,\bigl|\fsec(\eta_{0,i})-2\fsec(\eta_{1,i})\bigr|\\
&\quad\;
+\frac{h_{i+1}^{2}}{12}\,\bigl|\fsec(\eta_{1,i+1})-2\fsec(\eta_{0,i+1})\bigr|,\label{eq:TriangleIneq}
\end{align}
where $\Delta y_{i}^{-} = \bigl(\fsec(\eta_{0,i})-2\fsec(\eta_{1,i})\bigr) h_{i}^{2}/12$ and 
$\Delta y_{i}^{+} = \bigl(\fsec(\eta_{1,i+1})-2\fsec(\eta_{0,i+1})\bigr) h_{i+1}^{2}/12$
are the offsets with respect to $\f(x_{i})$ of the optimized segments~\eqref{eq:DeltaYMinimizers} at the left and right
of $x=x_{i}$, respectively; $\eta_{0,j}$ and $\eta_{1,j}$ lie in the
interval~$I_{j}$. Since we are dealing with functions twice
continuously differentiable in a closed interval, the absolute value
terms in~\eqref{eq:TriangleIneq} are bounded, according to the extreme
value theorem; therefore, if $h_{i}$ and $h_{i+1}$ decrease (finer partition $\part$), 
the discontinuity jumps at the knots of the partition also decrease. 
In the limit, as $\sup h_{i}\rightarrow 0$, $|\rightDeltaY i^{-}-\leftDeltaY i^{+}|\rightarrow 0$, i.e., continuity is satisfied.
Therefore the union of the independently-optimized segments $\mathrm{line}_{i}\rightarrow\ortproj\f$, 
which is the unique piecewise linear function satisfying both continuity ($\in \vecspace$) and minimization of the $\Ltwo$ error~\eqref{eq:phf-is-the-best}.
Consequently, if $N$ is large we may approximate
$\|\f-\ortproj\f\|_{\Ltwo(I_{i})}^{2}\approx\min\|\f-\mathrm{line}_{i}\|_{\Ltwo(I_{i})}^{2}$;
moreover, if $I_{i}$ is sufficiently small so that $\fsec$ is
approximately constant within it, $\fsec_{I_{i}}$, then we use Corollary~\ref{thm:MinErrSmallLineSegment} to get~\eqref{eq:approximate-best-ortproj}.

\ifCLASSOPTIONcaptionsoff
  \newpage
\fi



\balance 
%
%
%

\end{document}